\documentclass[11pt]{amsart}

\usepackage{epsfig, subfig, amsmath, amssymb, latexsym, amsfonts, xcolor, easybmat, bbm}

\newcommand\ol{\ensuremath{\overline}}

\newcommand\eop{{{\hfil \ensuremath \Box}}}

\newcommand\eps{\ensuremath {\varepsilon}}

\newenvironment{cor}{\subsection{}{\textbf {Corollary.}}\em}{}
\newenvironment{defn}{\subsection{}{\textbf {Definition.}}\em}{\smallskip}
\newenvironment{eg}{\subsection{}{\textbf {Example.}}}{\smallskip}

\newenvironment{lem}{\subsection{}{\textbf {Lemma.}}\em}{\smallskip}

\newenvironment{prop}{\subsection{}{\textbf {Proposition.}}\em}{\smallskip}

\newenvironment{rem}{\subsection{}{\textbf {Remark.}}}{\smallskip}
\newenvironment{thm}{\subsection{}{\textbf {Theorem.}}\em}{\smallskip}

\newenvironment{pf}{\noindent{\textbf {Proof.}}} {\begin{flushright}\eop \end{flushright}\smallskip}

\newcommand\fC{\ensuremath{\mathfrak C}}

\newcommand\fE{\ensuremath{\mathfrak E}}

\newcommand\fM{\ensuremath{\mathfrak M}}

\newcommand\fT{\ensuremath{\mathfrak T}}

\newcommand\cA{\ensuremath{\mathcal A}}
\newcommand\cB{\ensuremath{\mathcal B}}

\newcommand\cH{\ensuremath{\mathcal H}}

\newcommand\cJ{\ensuremath{\mathcal J}}
\newcommand\cK{\ensuremath{\mathcal K}}
\newcommand\cL{\ensuremath{\mathcal L}}
\newcommand\cM{\ensuremath{\mathcal M}}
\newcommand\cN{\ensuremath{\mathcal N}}

\newcommand\cS{\ensuremath{\mathcal S}}

\newcommand\cU{\ensuremath{\mathcal U}}

\newcommand\bbC{\ensuremath{\mathbb C}}
\newcommand\bbD{\ensuremath{\mathbb D}}

\newcommand\bbM{\ensuremath{\mathbb M}}
\newcommand\bbN{\ensuremath{\mathbb N}}

\newcommand\bbT{\ensuremath{\mathbb T}}

\newcommand\bbZ{\ensuremath{\mathbb Z}}

\newcommand\ttt{\ensuremath{\textsc}}

\newcommand\bofh{\ensuremath{\cB ( \cH)}}
\newcommand\kofh{\ensuremath{\cK ( \cH)}}

\newcommand\hilb{\ensuremath{\mathcal H}}
\newcommand\clos{\ensuremath{\textsc{clos}}\, }

\newcommand\norm{\ensuremath {\Vert}}

\newcommand\ran{\ensuremath {\mathrm {ran}}}

\newcommand\nul{\ensuremath {\mathrm {nul}}}

\newcommand\rank{\ensuremath{\mathrm{rank}\, }}

\newcounter{asst}
\newcounter{lab}
\newcounter{asstAA}
\newcounter{asstBB}
\newcounter{asstCC}
\newcounter{asstDD}
\newcounter{asstEE}
\newcounter{asstFF}
\newcounter{asstGG}
\newcounter{asstHH}
\newcounter{asstII}
\newcounter{asstJJ}
\newcounter{asstKK}
\newcounter{asstLL}
\newcounter{asstMM}
\newcounter{suppAA}
\definecolor{myred}{rgb}{0.6,0,0}
\definecolor{myblue}{rgb}{0,0.2,0.4}
\definecolor{mygreen}{rgb}{0.2,0.6, 0.5}

\newcommand\cniltwo{\ensuremath{\mathfrak{c}(\ttt{nil}_2)}}

\newcommand\niltwo{\ensuremath{\ttt{nil}_2}}

\newcommand\sqre{\ensuremath{(\ttt{SQ})}}

\begin{document}

%%%%%%%%%%%%%%%%%%%%%
%Title
%%%%%%%%%%%%%%%%%%%%%

\title{On commutators of square-zero Hilbert space operators}

%%%%%%%%%%%%%%%%%%%%%
% Authorship
%%%%%%%%%%%%%%%%%%%%%

\thanks{${}^1$ Research supported in part by NSERC (Canada)}
\thanks{${}^2$ Research supported in part by National Natural Science Foundation of China (No.: 12071174)}

\thanks{{\ifcase\month\or Jan.\or Feb.\or March\or April\or May\or
June\or
July\or Aug.\or Sept.\or Oct.\or Nov.\or Dec.\fi\space \number\day,
\number\year}}
\author
	[L.W. Marcoux]{{Laurent W.~Marcoux${}^1$}}
\address
	{Department of Pure Mathematics\\
	University of Waterloo\\
	Waterloo, Ontario \\
	Canada  \ \ \ N2L 3G1}
\email{LWMarcoux@math.uwaterloo.ca}

%%%%%%%%

\author
	[H. Radjavi]{{Heydar Radjavi}}
\address
	{Department of Pure Mathematics\\
	University of Waterloo\\
	Waterloo, Ontario \\
	Canada  \ \ \ N2L 3G1}
\email{hradjavi@uwaterloo.ca}

%%%%%%%%

\author
	[Y.H.~Zhang]{{Yuanhang~Zhang${}^2$}}
\address
	{School of Mathematics\\
	Jilin University\\
	Changchun 130012\\
	P.R. CHINA}
\email{zhangyuanhang@jlu.edu.cn}

%%%%%%%%%%%%%%%%%%%%%
% Abstract
%%%%%%%%%%%%%%%%%%%%%

\begin{abstract}
Let $\hilb$ be a complex, separable Hilbert space, and set $\cniltwo = \{ MN - NM : N, M \in \bofh, M^2 = 0 = N^2 \}$.   When $\dim\, \hilb$ is finite, we characterise the set  $\cniltwo$ and its norm-closure $\ttt{clos}(\cniltwo)$.  In the infinite-dimensional setting, we characterise the intersection of $\ttt{clos}(\cniltwo)$ with the set of biquasitriangular operators, and we exhibit an index obstruction to belonging to $\ttt{clos}(\cniltwo)$.

\end{abstract}

%%%%%%%%%%%%%%%%%%%%%
% Keywords and Subject Classification
%%%%%%%%%%%%%%%%%%%%%


\maketitle
\markboth{\textsc{  }}{\textsc{}}

%%%%%%%%%%%%%%%%%%%%%%%%%%%%%%%%%%%%%%%%%%
%%%%%%%%%%%%%%%%%%%%%%%%%%%%%%%%%%%%%%%%%%
% SECTION ONE
%%%%%%%%%%%%%%%%%%%%%%%%%%%%%%%%%%%%%%%%%%
%%%%%%%%%%%%%%%%%%%%%%%%%%%%%%%%%%%%%%%%%%

\section{Introduction} \label{section1}

%%%%%%%%

\subsection{} \label{sec1.01}
Throughout this paper, $\hilb$ will denote a complex, separable Hilbert space.  Denote by $\bofh$ the algebra of bounded linear operators acting on $\hilb$, and by $\kofh$ the closed-two sided ideal of $\bofh$ consisting of compact operators.   If $A \in \bofh$ (and more generally, if $A \in \cA$, a unital Banach algebra), we denote the spectrum of $A$ by $\sigma(A) := \{ \alpha \in \bbC: A - \alpha I \text{ is not invertible}\}$.   If $\sigma(A) = \{ 0\}$, we say that $A$ is \textbf{quasinilpotent}.

 Given $A, B \in \bofh$, the \textbf{commutator} of $A$ and $B$ is the element $[A, B] := A B - B A$.   When $\dim\, \hilb = n < \infty$, we find that $\hilb$ is isomorphic to $\bbC^n$, and a classical result of Albert and Muckenhoupt~\cite{AlbertMuckenhoupt1957} asserts that $T \in \bofh \simeq \bbM_n(\bbC)$ is a commutator if and only if the trace of $T$ is zero.    When $\dim\, \hilb = \infty$, a result of Brown and Pearcy~\cite{BrownPearcy1965}  shows that an operator $T \in \bofh$ is a commutator if and only if $T$ is not of the form $\lambda I + K$, where $0 \ne \lambda \in \bbC$ and $K \in \kofh$.

There is a significant literature which seeks to classify the set of commutators in various ideals of $\bofh$~\cite{Anderson1977, BeltitaPatnaikWeiss2014, Weiss1980, Weiss1986, Weiss2005}, as well as commutators (and their spans) in various operator algebras~\cite{DykemaKrishnaswamyUsha2018, DykemaSkripka2012, Fack1982, KaftalNgZhang2014, Marcoux2006, Marcoux2010, Robert2015}.  A fascinating but difficult open question is to completely characterise the set $[\kofh, \kofh] := \{ K L - L K: K, L \in \bofh \text{ compact}\}$.

In recent years there has been growing interest in characterising
\[
\mathfrak{c}(\Omega) := \{ [A, B] : A, B \in \Omega\}, \]
where $\Omega \subseteq \bofh$ is a fixed set of operators (e.g.~\cite{DrnovsekRadjaviRosenthal2002, RadjaviRosenthal2002, MarcouxRadjaviZhang2023, MarcouxRadjaviZhang2024}).   In the infinite-dimensional setting, this problem can often be intractable.  Part of the reason for this is that when $\Omega$ is \textbf{invariant under similarity} (i.e. $A \in \Omega$ implies that $S^{-1} A S \in \Omega$ for all invertible operators $S \in \bofh$), then so is $\mathfrak{c}(\Omega)$.   The problem is that  similarity-invariant subsets of $\bofh$ are in general poorly understood.   On the other hand, the work of Apostol, Fialkow, Herrero and Voiculescu~\cite{ApostolFialkowHerreroVoiculescu1984}  and others establishes an interesting framework for characterising the \emph{norm-closures} of similarity-invariant sets in terms of readily verifiable data (such as spectrum, rank, and semi-Fredholm index as defined below).    For this reason, when considering infinite-dimensional Hilbert spaces, we shall focus our attention on characterising the set $\ttt{clos} (\mathfrak{c}(\Omega))$.

For example, in the papers~\cite{MarcouxRadjaviZhang2023, MarcouxRadjaviZhang2024}, the authors studied the case where $\Omega = \fE := \{ E \in \bofh: E = E^2\}$, characterising $\mathfrak{c}(\fE)$ and $\ttt{clos}(\mathfrak{c}(\fE))$ when $\dim\, \hilb < \infty$, as well as the intersection of $\mathfrak{c}(\fE)$ with the set $\ttt{(bqt)}$ of biquasitriangular operators (defined in Paragraph~\ref{sec1.07} below) in the infinite-dimensional setting.

In the present paper, we consider a  second instance of this problem.   Let
\[
\niltwo := \{ M \in \bofh: M^2 = 0\} \]
denote the set of square-zero operators in $\bofh$.    When $\dim\, \hilb < \infty$, we are able to completely characterise both $\cniltwo$ and $\ttt{clos}(\cniltwo)$.   Paralleling our work on commutators of idempotents, when $\dim\, \hilb = \infty$, we are able to characterise the intersection of $\ttt{clos}(\cniltwo)$ with the set $\ttt{(bqt)}$, obtaining the surprising result that
\[
\ttt{clos}(\cniltwo) \cap \ttt{(bqt)}   = \ttt{clos}(\mathfrak{c}(\fE)) \cap \ttt{(bqt)}. \]

We next gather a number of definitions and relevant results that we shall require below.

%%%%%%%%%%%%%%%%%%%%%%%%%%%%%%%%%%%%%%%%%%

\begin{rem} \label{rem1.02}
We have obliquely alluded to the fact that the set $\niltwo$ is similarity-invariant.   Indeed, let $\hilb$ be a Hilbert space, and suppose that  $T \in \bofh$ lies in $\cniltwo$.   Choose $M, N \in \niltwo$ such that $T = [M, N]$.    If $S \in \bofh$ is invertible, then
\[
S^{-1} T S = S^{-1} [M, N] S = [S^{-1} M S, S^{-1} N S ], \]
and $(S^{-1} M S)^2 = 0 = (S^{-1} N S)^2$.  In other words, the following are equivalent.
\begin{itemize}
	\item{}
	$T \in \cniltwo$;
	\item{}
	there exists $S \in \bofh$ invertible such that $S^{-1} T S \in \cniltwo$;
	\item{}
	$\cS(T) := \{ R^{-1} T R : R \in \bofh \text{ invertible}\}  \subseteq \cniltwo$.
\end{itemize}	

Given $A, B \in \bofh$, let us write $A \sim B$ to indicate that $A$ and $B$ are similar.
Suppose that $\dim\, \hilb = n <\infty$.   Given $T \in \cB(\bbC^n)$, by considering the Jordan form of $T$, we see that we may write $T \sim B \oplus Q$, where $B$ is invertible and $Q$ is nilpotent.    The above observation implies that $T \in \cniltwo$ if and only if $B \oplus Q \in \cniltwo$, and our arguments below will often start with the assumption that $T$ is of the form $B \oplus Q$.

\smallskip

We point out that the same reasoning may be applied to the class $\sqre := \{ Z^2: Z \in \bofh\}$.   Thus the following are equivalent:
\begin{itemize}
	\item{}
	$T \in \sqre$;
	\item{}
	there exists $S \in \bofh$ invertible such that $S^{-1} T S \in \sqre$;
	\item{}
	$\cS(T) \subseteq \sqre$.
\end{itemize}	

A second well-studied equivalence relation on $\bofh$ is that of approximate unitary equivalence.   Two operators $A, B \in \bofh$ are said to be \textbf{approximately unitarily equivalent} if there exists a sequence $(U_n)_n$ of unitary operators acting on $\hilb$ such that $B = \lim_n U_n^* A U_n$.   As we shall see below (see Example~\ref{eg4.06}), although $\niltwo$ is invariant under approximate unitary equivalence, $\cniltwo$ is not.
\end{rem}

%%%%%%%%%%%%%%%%%%%%%%%%%%%%%%%%%%%%%%%%%%

\subsection{} \label{sec1.03}
In light of the comments above, in the finite-dimensional setting, nilpotent Jordan cells will play a central role in obtaining our characterisation of $\cniltwo$.   We denote by $J_n$ the usual Jordan cell in $\cB(\bbC^n)$, and more generally we shall write $J$ to denote a (finite) direct sum of Jordan cells.   At times we shall express an operator $T \in \cniltwo$ as a direct sum $T = T_1 \oplus T_2$ operators, each of which is a commutator of square-zero operators (acting on subspaces of the space $\hilb$ upon which $T$ acts).   In order to keep the notation as simple as possible, we shall say that $T_1, T_2$ and $T \in \cniltwo$, although the three operators may be acting on spaces of different dimensions.

An easy but useful example of an element in $\cniltwo$ is an operator of the form
\[
\begin{bmatrix} X & 0 \\ 0 & -X \end{bmatrix} = \left[ \begin{bmatrix} 0 & X \\ 0 & 0 \end{bmatrix}, \begin{bmatrix} 0 & 0 \\ I & 0 \end{bmatrix} \right], \]
where $X \in \bofh$ is arbitrary.   It is an easy exercise to verify that a finite direct sum of elements of $\cniltwo$ again lies in $\cniltwo$, and we shall implicitly apply this observation below.  In particular, this shows that $X \oplus -X \oplus 0 \in \cniltwo$ for all $X \in \bofh$.

%%%%%%%%%%%%%%%%%%%%%%%%%%%%%%%%%%%%%%%%%%

\begin{rem} \label{rem1.04}
A second simple but useful observation which we shall frequently use is the following.   If $T \in \cniltwo$, say $T = [M, N]$ where $M^2 = 0 = N^2$, then
\[
M T = - T M \text{     and      } N T = - T N. \]
\end{rem}

%%%%%%%%%%%%%%%%%%%%%%%%%%%%%%%%%%%%%%%%%%

\subsection{} \label{sec1.05}
Let $\hilb$ and $\cK$ be Hilbert spaces, $A \in \bofh$ and $B \in \cB(\cK)$.    We define the \textbf{Rosenblum operator}
\[
\begin{array}{rccc}
	\tau_{A, B}: & \cB(\cK, \hilb) & \to &\cB(\cK, \hilb) \\
		& X & \mapsto & A X - X B.
\end{array} \]

%%%%%%%%%%%%%%%%%%%%%%%%%%%%%%%%%%%%%%%%%%
		
The following result, due to Rosenblum, will prove very useful.  Its proof may be found in~\cite[Corollary 3.20]{Herrero1989}.

\begin{thm} \label{thm1.06} \emph{\textbf{(Rosenblum)}}

Let $\hilb$ and $\cK$ be Hilbert spaces, $A \in \bofh$ and $B \in \cB(\cK)$.   Then
\[
\sigma(\tau_{A, B}) = \sigma(A) - \sigma(B) := \{ \alpha - \beta: \alpha \in \sigma(A), \beta \in \sigma(B) \}. \]
In particular, if $\sigma(A) \cap \sigma(B) = \varnothing$, then $\tau_{A, B}$ is invertible.

Consequently, if $A \in \bofh$, $B \in \cB(\cK)$ and $Y \in \cB(\cK, \hilb)$, and if $\sigma(A) \cap \sigma(B) = \varnothing$, then
\[
\begin{bmatrix} A & Y \\ 0 & B \end{bmatrix} \text{      and      } \begin{bmatrix} A & 0 \\ 0 & B \end{bmatrix} \]
acting on $\hilb \oplus \cK$ are similar.
\end{thm}		

%%%%%%%%%%%%%%%%%%%%%%%%%%%%%%%%%%%%%%%%%%

\bigskip

Let us now define the set of biquasitriangular operators mentioned above, as well as the set of balanced operators which will play a role in Section~\ref{section4} of the paper.

%%%%%%%%%%%%%%%%%%%%%%%%%%%%%%%%%%%%%%%%%%

\subsection{} \label{sec1.07}
Let $\pi: \bofh \to \bofh/\kofh$ denote the canonical quotient map from $\bofh$ into the \textbf{Calkin algebra}.   Given $T \in \bofh$, the \textbf{essential spectrum} of $T$ is the spectrum of $\pi(T)$ in the Calkin algebra.    We say that $T$ is \textbf{semi-Fredholm} if $\pi(T)$ is either left- or right-invertible in the Calkin algebra, and that $T$ is \textbf{Fredholm} if $\pi(T)$ is invertible in $\bofh/\kofh$.  The \textbf{semi-Fredholm domain} of $T$ is the set
\[
\rho_{\text{s-F}}(T) := \{ \alpha \in \bbC: T - \alpha I \text{ is semi-Fredholm} \}, \]
and if $\alpha \in \rho_{\text{s-F}}(T)$, we define the \textbf{semi-Fredholm index} of $T-\alpha I$ as
\[
\ttt{ind} (T-\alpha I) := \ttt{nul} \, (T-\alpha I) - \ttt{nul}\, (T-\alpha I)^*. \]
It is well known that this is well-defined and lies in $\bbZ \cup \{-\infty, \infty\}$, as at least one of $\ttt{nul}\, (T-\alpha I)$ and $\ttt{nul}\, (T-\alpha I)^*$ is finite when $\alpha \in \rho_{\text{s-F}}(T)$.

An operator $T \in \bofh$ is said to be \textbf{biquasitriangular} if $\ttt{ind}\, (T-\alpha I) = 0$ whenever $\alpha \in \rho_{\text{s-F}}(T)$.   Although this is not the original definition of biquasitriangularity (see, e.g.~\cite{Halmos1970}), it is known to be equivalent to it~\cite[Corollary~5.5]{ApostolFoiasVoiculescu1974.02}, and it is the most practical form of the property with which to work.  We denote by $\ttt{(bqt)}$ the set of biquasitriangular operators in $\bofh$.

%%%%%%%%%%%%%%%%%%%%%%%%%%%%%%%%%%%%%%%%%%

The next definition played a crucial role in the study of $\mathfrak{c}(\fE)$ and its closure~\cite{MarcouxRadjaviZhang2024}, and we shall also require it here.

\bigskip

\begin{defn} \label{defn1.08}
Let $T \in \bofh$.   We say that $T$ is \textbf{balanced} if
\begin{enumerate}
	\item[(a)]
	$\sigma(T) = \sigma(-T)$;
	\item[(b)]
	whenever $\Omega_1, \Omega_2 \subseteq \bbC$ are disjoint open sets such that $\sigma(T) \subseteq \Omega_1 \cup \Omega_2$, then
	\[
	\dim\, \hilb(\Omega_1; T) = \dim\, \hilb(-\Omega_1; T), \]
	where $\hilb(\Omega_1;T)$ is the generalised eigenspace (i.e. the range of the corresponding Riesz idempotent $E(\Omega_1;T)) $corresponding to $\sigma(T) \cap \Omega_1$;
	\item[(c)]
	if $\alpha \in \bbC$, then $T - \alpha I$ is semi-Fredholm if and only if $T + \alpha I$ is semi-Fredholm, in which case
	\[
	\ttt{ind} (T - \alpha I) = \ttt{ind} (T + \alpha I). \]
\end{enumerate}
We denote the set of balanced operators on $\hilb$ by $\ttt{Bal}(\hilb)$.

We shall refer to an operator satisfying conditions (a) and (b) as \textbf{weakly balanced}, and we denote that set of operators by $\ttt{WBal}(\hilb)$.
\end{defn}

%%%%%%%%%%%%%%%%%%%%%%%%%%%%%%%%%%%%%%%%%%

\subsection{} \label{sec1.09}

We recall from Proposition~3.4 of~\cite{MarcouxRadjaviZhang2024} that if $\dim\, \hilb < \infty$, then
\[
\ttt{Bal}(\hilb) = \{ T \in \bofh: \alpha \in \sigma(T) \text{ implies that } -\alpha \in \sigma(T) \text{ and } \mu(\alpha) = \mu(-\alpha) \}, \]
where $\mu(\alpha)$ denotes the algebraic multiplicity of the eigenvalue $\alpha \in \sigma(T)$.

%%%%%%%%%%%%%%%%%%%%%%%%%%%%%%%%%%%%%%%%%%

%%%%%%%%%%%%%%%%%%%%%%%%%%%%%%%%%%%%%%%%%%
%%%%%%%%%%%%%%%%%%%%%%%%%%%%%%%%%%%%%%%%%%
% SECTION TWO
%%%%%%%%%%%%%%%%%%%%%%%%%%%%%%%%%%%%%%%%%%
%%%%%%%%%%%%%%%%%%%%%%%%%%%%%%%%%%%%%%%%%%

\vskip 1.5 cm

\section{Structure results for elements of $\cniltwo$} \label{section2}

%%%%%%%%

%%%%%%%%%%%%%%%%%%%%%%%%%%%%%%%%%%%%%%%%%%

\subsection{} \label{sec2.01}
In this section, we establish a number of  results concerning elements of $\cniltwo$.     In particular, we shall characterise those invertible operators and those normal operators which lie in $\cniltwo$.

Our first result will be of particular importance in the finite-dimensional setting where,  by considering Jordan decompositions, we see that every operator in $\cB(\bbC^n)$ is similar to a direct sum of an invertible operator $B$ with a nilpotent operator $J$.

%%%%%%%%%%%%%%%%%%%%%%%%%%%%%%%%%%%%%%%%%%

\begin{prop} \label{prop2.02}
Let $\hilb$ be a Hilbert space  and suppose that $T \sim A \oplus B \in \bofh$, where $\sigma(A^2) \cap \sigma(B^2) = \varnothing$.  The following statements are equivalent:
\begin{enumerate}
	\item[(a)]
	$T \in \cniltwo$; and
	\item[(b)]
	$A \in \cniltwo$ and $B \in \cniltwo$.
\end{enumerate}
\end{prop}	

\begin{pf}

By Remark~\ref{rem1.02}, we may suppose without loss of generality that $T = A \oplus B$.
\begin{enumerate}
	\item[(a)] implies (b). \ \ \
	Suppose that $\hilb = \hilb_A \oplus \hilb_B$, that $A \in \cB(\hilb_A)$, $B \in \cB(\hilb_B)$, and $\sigma(A^2) \cap \sigma(B^2) = \varnothing$.
	Suppose furthermore that $T = (A \oplus B) \in \cniltwo$, and choose square-zero operators $M,N \in \bofh$ such that $T=MN-NM$.
	Routine calculations show that
	\[
	T^2 M = M N M N M = M T^2, \]
	and
	\[
	T^2 N = N M N M N = N T^2.\]	
	With respect to the decomposition $\hilb = \hilb_A \oplus \hilb_B$, we write
	\[
	M=\begin{bmatrix}M_1& M_2 \\ M_3& M_4 \end{bmatrix},~~N=\begin{bmatrix}N_1& N_2 \\ N_3& N_4 \end{bmatrix}.\]
	By comparing the $(1,2)$ and the $(2,1)$ entries  of $T^2M$ and $MT^2$, we conclude that
	\[
	A^2 M_2=M_2 B^2,~~B^2 M_3= M_3 A^2.\]
	Since $\sigma(A^2) \cap \sigma(B^2) =\varnothing$, by Theorem~\ref{thm1.06}, the corresponding Rosenblum operators
	$\tau_{A^2, B^2}$ and $\tau_{B^2, A^2}$ are injective.  From this we may conclude that  $M_2 = 0 = M_3$, so that
	$M = M_1 \oplus M_4$.  Similarly $N = N_1 \oplus N_4$.  Since $M^2 = 0 = N^2$, we conclude that
	$M_1, M_4, N_1$ and $N_4$ are all nilpotent of order at most two.
	
	But $T = [M, N]$, and thus $A = [M_1, N_1]$ and $B = [M_4,  N_4]$ lie in $\cniltwo$.
	\smallskip
	
	\item[(b)] implies (a). \ \ \
	If $A = [M_1, N_1]$ and $B = [M_2, N_2]$ where $M_i^2 = 0 = N_i^2$, $i = 1, 2$, then $(M_1 \oplus M_2)^2 = 0 = (N_1 \oplus N_2)^2$ and
	\[
	T = [(M_1 \oplus M_2), (N_1 \oplus N_2)]. \]
\end{enumerate}
\end{pf}	

%%%%%%%%%%%%%%%%%%%%%%%%%%%%%%%%%%%%%%%%%%

\begin{prop} \label{prop2.03}
Let $\hilb$ be a Hilbert space and suppose that $B \in \bofh$ is invertible.  The following conditions are equivalent.
\begin{enumerate}
	\item[(a)]
	$B \in \cniltwo$.
	\item[(b)]
	There exists an invertible operator $A_0$ such that $B \sim A_0 \oplus -A_0$.
\end{enumerate}		
\end{prop}

\begin{pf}
\begin{enumerate}
	\item[(a)]
	Suppose that $B \in \cniltwo$.  Let $B = [M, N]$ where $M, N \in \bofh$ satisfy $M^2 = 0 = N^2$.

	Decompose $\hilb = \hilb_1 \oplus \hilb_2$, where $\hilb_1 = \ker\, M$. Relative to this decomposition, we write
	\[
	M=\begin{bmatrix} 0 & A \\ 0 & 0 \end{bmatrix}, \ \ \ \ \ N = \begin{bmatrix} W & X \\ Y & Z \end{bmatrix}. \]
	It follows that
	\[
	B = [M, N] = \begin{bmatrix} A Y & A Z - W A \\ 0 & - Y A \end{bmatrix}. \]

	\smallskip

	If $\dim\, \hilb < \infty$, then the invertibility of $B$ implies that of both $A Y$ and $-Y A$, which in turn implies (by considering the
	ranks of $AY$ and of $-YA$) that $\dim\, \hilb_1 = \dim\, \hilb_2$ and thus $\dim\, \hilb$ is even!
	
	If $\dim\, \hilb = \infty$, then the invertibility of $B$ implies that $M$ has infinite rank.  Since $M^2 = 0$, it clearly has infinite nullity,
	and thus $\dim\, \hilb_1 = \dim\, \hilb_2$. (The same argument shows that the rank and nullity of $N$ are also infinite.)

	In either the finite-dimensional or the infinite-dimensional setting, it follows from the invertibility of $B$ that $A Y$ is left-invertible
	and that $Y A$ is right-invertible.   In particular, $\ker\, AY = \{ 0\}$ and $\ran\, Y A = \hilb_2$.  From this it follows that
	$\ker\, Y = \{ 0\}$ and $\ran\, Y = \hilb_2$, so that $Y$ is a (continuous) bijection.   By the Open Mapping Theorem, $Y$ is invertible.

	Let $S := I \oplus Y$ and consider $M_0 := S^{-1} M S$, $N_0 := S^{-1} N S$.   Then, for appropriate choices of
	$A_0, W_0, X_0$ and $Z_0$, $M_0$ and $N_0$ are of the form
	\[
	M_0 = \begin{bmatrix} 0 & A_0 \\ 0 & 0 \end{bmatrix}, \ \ \ \ \ N_0 = \begin{bmatrix} W_0 & X_0 \\ I  & Z_0 \end{bmatrix}. \]
	Since $N^2 = 0$, we see that $Z_0 = -W_0$ and that $X_0 = -Z_0^2$.   Thus we may write
	\[
	N_0 = \begin{bmatrix} W_0 & - W_0^2 \\ I & - W_0 \end{bmatrix}. \]
	Next, note that
	\[
		S^{-1}  B  S =  S^{-1} [M, N] S = [M_0, N_0] = \begin{bmatrix} A_0 & -A_0 W_0 - W_0 A_0 \\ 0 &  -A_0 \end{bmatrix}. \]
	Set $R = \begin{bmatrix} W_0 & I \\ I & 0 \end{bmatrix}$, so that $R$ is invertible with inverse $R^{-1} = \begin{bmatrix} 0 & I \\ I & -W_0 \end{bmatrix}$.
	A routine calculation shows that
	\[
	(SR)^{-1} B (SR) = R^{-1} (S^{-1} B S) R =  (-A_0) \oplus A_0. \]	

	The last statement is clear.
	\item[(b)]
	This was demonstrated in Paragraph~\ref{sec1.03}.
\end{enumerate}	
\end{pf}

%%%%%%%%%%%%%%%%%%%%%%%%%%%%%%%%%%%%%%%%%%

\begin{defn} \label{defn2.03.01}
Let $T\in \bofh$.  We say that $T$ is \textbf{strongly irreducible} if there is no nontrivial idempotent in the commutant $\{ T\}^\prime  := \{ X \in \bofh: T X = X T\}$ of $T$.
\end{defn}

\smallskip

The reader is referred to the monograph of Jiang and Wang~\cite{JiangWang1998} for more information about strongly irreducible operators.

\bigskip

%%%%%%%%%%%%%%%%%%%%%%%%%%%%%%%%%%%%%%%%%%

\begin{cor} \label{cor2.03.02}
If $T \in \bofh$ is an invertible strongly irreducible operator, then $T\notin \cniltwo$.
\end{cor}

\begin{pf}
This follows immediately from Proposition~\ref{prop2.03}.
\end{pf}

%%%%%%%%%%%%%%%%%%%%%%%%%%%%%%%%%%%%%%%%%%

The following proposition is derived from the proof of Corollary 2 of Gellar's paper~\cite{Gellar1969} and Corollary~\ref{cor2.03.02}.  Note that every weighted shift $W$ acting on a Hilbert space is unitarily equivalent to $-W$.

\begin{prop} \label{prop2.03.03}
Let $W \in B(\cH)$ be an injective bilateral weighted shift.   If there exist $0 < R_1 < R_2$ such that $\sigma(W )=\{z: R_1\leq |z|\leq R_2\}$,
then $W$ is strongly irreducible.  As such, $W \not \in \cniltwo$.
\end{prop}

%%%%%%%%%%%%%%%%%%%%%%%%%%%%%%%%%%%%%%%%%%

\begin{eg} \label{eg2.03.04}
Let $\{ e_n\}_{n \in \bbZ}$ be an orthonormal basis for the Hilbert space $\hilb$, and consider the weighted shift $W$ with weight sequence $(w_n)_{n \in \bbZ}$, where $w_n = \begin{cases} 2 & \mbox{ if } n < 0 \\ 1 & \mbox{ if } n \ge 0 \end{cases}$.

Then $W$ is invertible, and the spectrum of $W$ is easily seen to be $\sigma(W)=\{z:1\leq |z|\leq 2\}$.   By Proposition~\ref{prop2.03.03}, $W$ is strongly irreducible and  thus $W\notin \cniltwo$, although $W$ is invertible and it is similar to $-W$.
\end{eg}

%%%%%%%%%%%%%%%%%%%%%%%%%%%%%%%%%%%%%%%%%%

\bigskip

Proposition~\ref{prop2.03} implies that if $B \in \cniltwo$ is invertible, then $B \sim -B$.    In fact, the condition $T \sim -T$ is equivalent to the condition $T \in \cniltwo$ when $T$ is normal, as we shall now demonstrate.

\begin{prop}\label{prop2.04}
Let $\hilb$ be a Hilbert space and suppose that $T \in \bofh$ is a normal operator.  If $T \sim -T$, then $T \in \cniltwo$.
\end{prop}

\begin{pf}	
First note that as $T$ and $-T$ are both normal, if $T$ is similar to $-T$, then $T$ is unitarily equivalent to $-T$ by~\cite[Problem~192]{Halmos1982}.

The next part of the argument is almost identical to that of \cite[Proposition~3]{DrnovsekRadjaviRosenthal2002}.
Let $E$ be the spectral measure of the normal operator $T$, and let
\[\Omega_1=\{z\in \mathbb{C}: z\neq 0, 0\leq \text{arg}z<\pi\},~~\Omega_2=\{z\in \mathbb{C}: z\neq 0, \pi\leq \textup{arg}z<2\pi\}.\]
Since $T$ is unitarily equivalent to $-T$, $T|_{E(\Omega_2)}$ is (unitarily) similar to $-T|_{E(\Omega_1)}$.
Thus, without loss of generality, one can assume that
$T$ has the form $A\oplus -A \oplus 0$, with respect to the decomposition $\cH=\textup{Ran }E(\Omega_1)\oplus \textup{Ran }E(\Omega_2)\oplus \textup{Ran }E({0})$.   As seen in Paragraph~\ref{sec1.03}, it easily follows that $T \in \cniltwo$.
\end{pf}

%%%%%%%%%%%%%%%%%%%%%%%%%%%%%%%%%%%%%%%%%%

\begin{lem} \label{lem2.05}
Let $\hilb$ be a Hilbert space and $T \in \bofh$ be a normal operator.   Suppose that $T = [M, N]$, where $M^2=0$, and $N^2=0$.  Relative to $\hilb = \ker\, M \oplus (\ker\, M)^\perp$, write
\[
M=\begin{bmatrix}0 & A\\ 0& 0\\ \end{bmatrix}, \ \ \ \ \ \ \ \ \ \ N=\begin{bmatrix}W & X\\ Y& Z\\ \end{bmatrix},\]
where $A$ is injective.   If $YA$ is injective, then $AZ-WA=0$.
\end{lem}

\begin{pf}
Note that $NM$ commutes with $T$. Since $T$ is normal, this implies that
\[NM=\begin{bmatrix}0 & WA\\ 0& YA\\ \end{bmatrix}~~\textup{and}~~T^*=\begin{bmatrix}(AY)^* & 0\\ (AZ-WA)^* & -(YA)^*\\ \end{bmatrix}\]
commute. In particular,
\[YA(AZ-WA)^*=0,\]
so that $(AZ-WA)^*=0$ by hypothesis.
\end{pf}

%%%%%%%%%%%%%%%%%%%%%%%%%%%%%%%%%%%%%%%%%%

\begin{thm} \label{thm2.06}
Let $\hilb$ be a Hilbert space and $T \in \bofh$ be normal.   The following conditions are equivalent.
\begin{enumerate}
	\item[(a)]
	$T$ is unitarily equivalent to $-T$; i.e. there exists a unitary operator $U \in \bofh$ such that $-T = U^* T U$.
	\item[(b)]
	$T \sim -T$.
	\item[(c)]
	$T \in \cniltwo$.
\end{enumerate}	
\end{thm}

\begin{pf}
\begin{enumerate}
	\item[(a)] implies (b).
	This is trivial, since each unitary operator is invertible.
	\item[(b)] implies (c).
	This is Proposition~\ref{prop2.04}
	\item[(c)] implies (a).
 	By Zorn's lemma there is a maximal subspace $\mathcal{L}$ invariant under both $M$ and $N$ on which the respective
	restrictions of $M$ and $N$ commute. (It is possible that $\cL = \{ 0\}$.) Relative to the decomposition $\hilb = \cL \oplus \cL^\perp$, $M$ and
	$N$ have the forms
	\[
	M= \begin{bmatrix}M_1 & *\\ 0& M_2\\ \end{bmatrix}
	 ,~~\textup{and}~~N=\begin{bmatrix}N_1 & *\\ 0& N_2\\ \end{bmatrix}.\]

	We may then apply the same argument to the adjoints of the nilpotents  $M_2$ and $N_2$.   This leads to a decomposition of $\cL^\perp$ as $\cL^\perp = \cM \oplus \cN$ where, with respect to the decomposition $\hilb = \cL \oplus \cM \oplus \cN$, we have:
	\[
	M=\begin{bmatrix} M_1&*&*\\
	0&M_0&*\\
	0&0&M_3 \\
	\end{bmatrix}~~\textup{and}~~
	N=\begin{bmatrix}
	N_1&*&*\\
	0&N_0&*\\
	0& 0& N_3\\
	\end{bmatrix}. \]
	Recall that the kernel and the range of a normal operator are reducing for that operator.   Since $T$ is normal and $T=MN-NM$, it follows that
	\[
	T=\begin{bmatrix}
	0&0&0\\
	0&T_0&0\\
	0& 0&0\\
	\end{bmatrix}.\]
	where  $T_0$  is the commutator of $M_0$ and $N_0$.  Note that $M_0$ and $N_0$ cannot have common kernel (by the maximality of $\cL$ and $\cN$).  Neither can their adjoints.

	We now write
	 \[
	 M_0=\begin{bmatrix}0 & A\\ 0& 0\\ \end{bmatrix},~~\textup{and}~~N_0=\begin{bmatrix}W & X\\ Y& Z\\ \end{bmatrix},\]
	relative to the decomposition $\mathcal{M}=\mathcal{M}_1\oplus (\mathcal{M}\ominus \mathcal{M}_1)$ with $\mathcal{M}_1=\ker M_0$.
	Then $A$ is injective and with respect to the decomposition $\cM = \cM_1 \oplus (\cM \ominus \cM_1)$, $T_0$ has the form
	\[
	\begin{bmatrix}AY&AZ-WA\\ 0&-YA\\ \end{bmatrix}.\]

	We claim that $Y$ is injective. Otherwise, $Y$ has nonzero kernel $\mathcal{M}_2$
	contained in $\mathcal{M}_1$. Note that $N_0^2=0$ implies $YW+ZY=0$, and thus
	$\mathcal{M}_2$ is invariant under $W$. Then with respect to the decomposition $\mathcal{M}=\mathcal{M}_2\oplus (\mathcal{M}_1\ominus
	\mathcal{M}_2)\oplus (\mathcal{M}\ominus \mathcal{M}_1)$, we have
	\[
	M_0=\begin{bmatrix} 0&0&A_1\\
	0&0&A_2\\
	0&0&0 \\
	\end{bmatrix}~~\textup{and}~~
	N_0=\begin{bmatrix}
	W_1&W_2&X_1\\
	0&W_3&X_2\\
	0& Y_1& Z\\
	\end{bmatrix}.\]
	Note that $\mathcal{L}\oplus \mathcal{M}_2$ is an invariant subspace under both $M$ and $N$ on which the respective
	restrictions of $M$ and $N$ commute.  This contradicts  the maximality of $\mathcal{L}$. Thus, $Y$ is injective.

	By Lemma \ref{lem2.05}, $AZ-WA=0$. Then
	\[T_0=\begin{bmatrix}AY&0\\ 0&-YA\\ \end{bmatrix}.\]
	Note that $Y(AY)=(YA)Y$, $A(YA)=(AY)A$, and $A,Y$ are injective, and thus by~\cite[Lemma 4.1]{ChenHerreroWu1992},
	$AY$ is unitarily equivalent to $YA$. Therefore, 	$T_0$ is unitarily equivalent to $-T_0$,
	and hence $T$ is unitarily equivalent to $-T$.
\end{enumerate}	
\end{pf}

%%%%%%%%%%%%%%%%%%%%%%%%%%%%%%%%%%%%%%%%%%

\begin{eg} \label{eg2.07}
By Example~\ref{eg2.03.04}, we saw that there exists an invertible operator $T \in \bofh$ which satisfies $T \sim -T$, and yet  $T$ does not lie in $\cniltwo$.

Let us now show that if the underlying Hilbert space $\hilb$ is infinite-dimensional, then there exists an operator
$T \in \cniltwo$ such that $T$ is not similar to $-T$.  The inspiration for this example comes from a matrix decomposition for products
of two square-zero operators found in the paper of Novak~\cite{Novak2008}.

\bigskip

Let $B \in \bofh$ be a positive, compact, injective, diagonal operator; say $B = \textsc{diag} (b_n)_n$.  (For example, one might choose $B=\textsc{diag} (\frac{1}{n})$.)

Consider $M = \begin{bmatrix} 0 & 0 & 0 \\ I & 0 & 0 \\ B & 0 & 0 \end{bmatrix}$ and $N= \begin{bmatrix} 0 & 0 & I \\ 0 & 0 & 0 \\ 0 & 0 & 0 \end{bmatrix} \in \cB(\hilb \oplus \hilb \oplus \hilb)$.   Clearly $M^2 = 0 = N^2$.

Let
\[
T := [M, N]  = \begin{bmatrix} - B & 0 & 0 \\ 0 & 0 & I \\ 0 & 0 & B \end{bmatrix}. \]

Suppose that $S = \begin{bmatrix} S_1 & S_2 & S_3 \\ S_4 & S_5 & S_6 \\ S_7 & S_8 & S_9 \end{bmatrix}$ satisfies $T S =  - S T$.   Then
\[
\begin{bmatrix}
	-B S_1 & -B S_2 & -B S_3  \\ S_7 & S_8 & S_9  \\ B S_7& B S_8 & B S_9 \end{bmatrix}
=
\begin{bmatrix}
	S_1 B & 0 & -(S_2 + S_3 B) \\ S_4 B & 0 & -(S_5 + S_6 B) \\ S_7 B & 0 & - (S_8 + S_9 B) \end{bmatrix}. \]

Since $B$ is injective, $- B S_2 = 0$ implies that $S_2 = 0$, and obviously $S_8 = 0$.

\smallskip

Next, note that
\begin{itemize}
	\item{}
	$B S_9 = - S_8 - S_9 B = -S_9 B$.   With respect to the basis that diagonalises $B$, write the corresponding matrix of $S_9$ as $[S_9] = [x_{ij}]$.   The equation $B S_9 = - S_9 B$ then implies that for all $1 \le i, j$, $b_i x_{ij} = - x_{ij} b_j$, or equivalently, $(b_i + b_j) x_{ij} = 0$.   Since $b_i + b_j > 0$ by hypothesis, we find that $x_{ij} = 0$ for all $1 \le i,j$.  That is, $S_9 = 0$.
	\item{}
	$S_7 = S_4 B$.
\end{itemize}
Thus
\[
S = \begin{bmatrix} S_1 & 0 & S_3 \\ S_4 & S_5 & S_6 \\ S_4 B & 0 & 0 \end{bmatrix}. \]

\bigskip

But $S_4 B$ acts on an infinite-dimensional space, and it is compact because $B$ is compact.    Thus $S$ is not invertible.    (If $S T = I$, then the third row of $S$ multiplied by the third column of $T$ is compact, contradicting the hypothesis that it should be $I$.)
\end{eg}

%%%%%%%%%%%%%%%%%%%%%%%%%%%%%%%%%%%%%%%%%%

\bigskip

Despite the fact that $T \in \cniltwo$ does not imply that $T \sim -T$, we now establish a number of useful results that illustrate a strong connection between $T$ and $-T$ when $T \in \cniltwo$.   The first such result establishes a link between $\cniltwo$ and $\ttt{Bal}(\hilb)$.

\begin{thm} \label{thm2.08}
Let $T \in \cniltwo$.   For all $0 \ne \alpha \in \bbC$ and $k \ge 1$,
\[
\nul \, (T-\alpha I)^k = \nul\, (T+\alpha I)^k. \]
\end{thm}

%%%%%%%%%%

\begin{pf}
$\bullet$ \ \ \textsc{Step one.}\ \ \

Let $T = [M, N] = MN - NM$, where $M^2 = 0 = N^2$.   Let $\hilb_1 := \ker\, M$ and $\hilb_2 = \hilb_1^\perp$,  so that relative to the decomposition $\hilb = \hilb_1 \oplus \hilb_2$, write $M = \begin{bmatrix} 0 & A \\ 0 & 0 \end{bmatrix}$ and $N = \begin{bmatrix} W & X \\ Y & Z \end{bmatrix}$.  Set $L := MN$ and $R := NM$, and observe that
\[
 	LR = (MN)(NM) = M N^2 M = 0 \]
and $RL = 0$ as well.

Now, $L = \begin{bmatrix} AY & AZ \\ 0 & 0 \end{bmatrix}$ and $R = \begin{bmatrix} 0 & WA \\ 0 & YA \end{bmatrix}$.  If $0 \ne \alpha \in \sigma_p(L)$, then $\alpha \in \sigma_p(AY)$.  Similarly, if $0 \ne \alpha \in \sigma_p(R)$, then $\alpha \in \sigma_p(YA)$.

%%%%%%%%%%%%%%%%%%%%%
\smallskip
%%%%%%%%%%%%%%%%%%%%%

$\bullet$\ \ \textsc{Step two.}\ \ \
Let $k \in \bbN$ and $0 \ne \alpha \in \bbC$.

\begin{enumerate}
	\item[(II.1)]
	Note  that if $0 \ne w = \begin{bmatrix} w_1 \\ w_2 \end{bmatrix} \in \ker\, (L - \alpha I)^k$, then  (for an appropriate operator $L_{k, \alpha, 2}$)
	\[
	\begin{bmatrix} 0 \\ 0 \end{bmatrix} = (L - \alpha I)^k w = \begin{bmatrix} (AY - \alpha I)^k w_1 + L_{k, \alpha, 2} w_2 \\ (-\alpha)^k \ w_2 \end{bmatrix}, \]
	implying that $w_2 = 0$.
	\item[(II.2)]
	Similarly, if $0 \ne w \in \ker\, (R+\alpha I)^k$, then (for an appropriate operator $R_{k, \alpha, 2}$)
	\[
	\begin{bmatrix} 0 \\ 0 \end{bmatrix} = (R+ \alpha I)^k w = \begin{bmatrix} \alpha^k  w_1 + R_{k, \alpha, 2}\  w_2 \\  (YA + \alpha)^k w_2 \end{bmatrix}. \]
	It follows that $w_2 \ne 0$, for otherwise $\alpha^k w_1 = 0$ and hence $w_1 = 0$, a contradiction.
\end{enumerate}	

%%%%%%%%%%%%%%%%%%%%%
\bigskip
%%%%%%%%%%%%%%%%%%%%%
	
$\bullet$\ \ \textsc{Step three.}\ \ \ 	

Now let $0 \ne \alpha \in \sigma_p(T)$, and let $z \in \ker\, (T-\alpha I)^k$.
Then (for an appropriate operator $T_{k, \alpha, 2}$)
\[
\begin{bmatrix} 0 \\ 0 \end{bmatrix}
	= (T - \alpha I)^k z
	= \begin{bmatrix} AY - \alpha I & A Z - W A \\ 0 & - YA - \alpha I \end{bmatrix}^k \ \begin{bmatrix} z_1 \\ z_2 \end{bmatrix}
	= \begin{bmatrix} (AY - \alpha I)^k z_1 + T_{k, \alpha, 2}\  z_2 \\  (-YA-\alpha I)^k z_2 \end{bmatrix}. \]
	
This leaves us with two possibilities.   Either
\begin{itemize}
	\item{}
	$z_2 \ne 0$, and thus $z_2 \in \ker\, (YA + \alpha I)^k$; or
	\item{}
	$z_2 = 0$, in which case $z_1 \ne 0$ and $(AY - \alpha I)^k z_1 = 0$.
\end{itemize}

%%%%%%%%%%%%%%%%%%%%%
\bigskip %\vfill\newpage
%%%%%%%%%%%%%%%%%%%%%
	
$\bullet$\ \ \textsc{Step four.}\ \ \ 	
Conversely, suppose that $0 \ne y_2 \in \hilb_2$ and that $(Y A + \alpha I)^k y_2 = 0$.

Let $y_1 := -(\alpha^{-k}) R_{k, \alpha, 2} \ y_2$, where $R_{k, \alpha, 2}$ is defined as above,  and  set $y = \begin{bmatrix} y_1 \\ y_2 \end{bmatrix}$.    Clearly $y \ne 0$ as $y_2 \ne 0$.  Then
\[
(R + \alpha I)^k \begin{bmatrix} y_1 \\ y_2 \end{bmatrix}
%	= \begin{bmatrix} 0 & W A \\ 0 & Y A  \end{bmatrix} \begin{bmatrix} y_1 \\ y_2 \end{bmatrix}
	=\begin{bmatrix} \alpha^k  y_1 + R_{k, \alpha, 2} \ y_2 \\ (YA + \alpha)^k  y_2 \end{bmatrix}
	=  \begin{bmatrix}  \alpha^k  y_1 + R_{k, \alpha, 2} \ y_2  \\ 0 \end{bmatrix}
	= \begin{bmatrix} 0 \\ 0 \end{bmatrix}. \]

That is, $(R + \alpha I)^k y = 0$.  Note that $(R + \alpha I)^k L y = L (R + \alpha I)^k  y = L 0 = 0$.   Hence $L y \in \ker \, (R+\alpha I)^k$.    But then
\[
L y = \begin{bmatrix} A Y & A Z \\ 0 & 0 \end{bmatrix} \begin{bmatrix} y_1 \\ y_2 \end{bmatrix} = \begin{bmatrix} (AY) y_1 + (AZ) y_2 \\ 0 \end{bmatrix}. \]
By (II.2) of \textsc{Step two} above, we must have $Ly = 0$.

\bigskip

A tedious but routine calculation shows that
\begin{align*}
(T-\alpha I)^k
	&= (L - (R + \alpha I))^k \\
	&= L^k -  \binom{k}{1} \alpha L^{k-1} + \binom{k}{2} \alpha^2 L^{k-2} + \cdots \\
	& \ \ \ \ \ \ \ \ \ \ \ \ \ \ \ \cdots + (-1)^{k-1} \binom{k}{k-1} \alpha^{k-1} L^{1} + (-1)^k (R+\alpha I)^k.
\end{align*}	
From this it follows that
\[
(T - \alpha I)^k y  = (L - (R + \alpha I))^k y =  0, \]
so that $y \in \ker \, (T - \alpha I)^k$.

More generally, if $y_2^{(1)}, y_2^{(2)}, \ldots, y_2^{(m)}$ is a linearly independent set of vectors in the kernel of $(YA + \alpha I)^k$, then with $y_1^{(j)} := (-\alpha)^{-k} R_{k, \alpha, 2} \ y_2^{(j)}$, $1 \le j \le m$, we see that
\[
y^{(j)} := \begin{bmatrix} y_1^{(j)} \\ y_2^{(j)} \end{bmatrix}, \ \ \ 1 \le j \le m \]
is a linearly independent family of vectors in $\ker \, (T - \alpha I)^k$.

%%%%%%%%%%%%%%%%%%%%%
\bigskip
%%%%%%%%%%%%%%%%%%%%%
	
$\bullet$\ \ \textsc{Step five.}\ \ \ 	Suppose that $0 \ne x_1 \in \hilb_1$ and that $(AY - \alpha I)^k x_1 = 0$.

Let $x := \begin{bmatrix} x_1 \\ 0 \end{bmatrix}$.   Then
\[
(L - \alpha I)^k x = 0. \]
% \begin{bmatrix} A Y & A Z \\ 0 & 0 \end{bmatrix} \begin{bmatrix} x_1 \\ 0 \end{bmatrix} = \begin{bmatrix} (AY) x_1 \\ 0 \end{bmatrix} = \alpha \begin{bmatrix} x_1 \\ 0 \end{bmatrix}. \]
A routine calculation similar to the one above shows that
\begin{align*}
(T-\alpha I)^k
	&= ( (L-\alpha I) - R)^k \\
	&=  (L-\alpha I)^ k   \\
	&\ \ \  +  (-1)^k \left ( \binom{k}{1} \alpha^{k-1} R + \binom{k}{2} \alpha^{k-2} R^2 + \cdots \cdots + \binom{k}{k-1} \alpha^1 R^{k-1} + R^k \right).
\end{align*}	

From this and the fact that $R x = 0$ it follows that
\[
(T - \alpha I)^k x = (L-\alpha I)^k  x + (-1)^k \sum_{j=1}^k \binom{k}{j} \alpha^{k-j} R^j x = \begin{bmatrix} 0 \\ 0 \end{bmatrix}. \]
Hence $x \in \ker\, (T-\alpha I)^k$.   More generally, if $x_1^{(1)}, x_1^{(2)}, \ldots, x_1^{(\ell)}$ is a linearly independent set of vectors in $\ker \, (AY - \alpha I)^k$, then
\[
x^{(j)} := \begin{bmatrix} x_1^{(j)} \\ 0 \end{bmatrix}, \ \ \ 1 \le j \le \ell \]
is a linearly independent family of vectors in $\ker\, (T-\alpha I)^k$.

%%%%%%%%%%%%%%%%%%%%%
\bigskip
%%%%%%%%%%%%%%%%%%%%%
	
$\bullet$\ \ \textsc{Step six.}\ \ \ 	
We also observe that $\{ x^{(1)}, x^{(2)}, \ldots, x^{(\ell)}, y^{(1)}, y^{(2)}, \ldots, y^{(m)} \}$  is linearly independent.   Indeed, if $\eta_i, \zeta_j \in \bbC$, $1 \le i \le \ell, 1 \le j \le m$ and
\[
\sum_{i=1}^\ell \eta_i x^{(i)} + \sum_{j=1}^m \zeta_j y^{(j)} = 0, \]
then by considering the second coordinate of these sums, we see that $\zeta_j = 0$, $1 \le j \le m$, and hence by considering the remaining first coordinates, $\zeta_i = 0$, $1 \le i \le \ell$.

\bigskip

By considering the case where  $\{ x^{(1)}, x^{(2)}, \ldots, x^{(\ell)}\}$ is a basis for $\ker\, (AY -\alpha I)^k$ and $\{ y^{(1)}, y^{(2)}, \ldots, y^{(m)} \}$ is a basis for $\ker\, (YA + \alpha I)^k$, we conclude that
\[
\nul\, (T-\alpha I)^k \ge \nul\, (AY - \alpha I)^k + \nul\, (YA + \alpha I)^k. \]
The reverse in equality follows easily from the form of $(T-\alpha I)^k$ exhibited in \textsc{Step three}.

%%%%%%%%%%%%%%%%%%%%%
\bigskip
%%%%%%%%%%%%%%%%%%%%%
	
$\bullet$\ \ \textsc{Step seven.}\ \ \ 	

Consider the following equation from  Horn and Johnson~\cite[Theorem1.3.20]{HornJohnson1990}:
\[
\begin{bmatrix} I & - A \\ 0 & I \end{bmatrix} \ \begin{bmatrix} A Y & 0 \\ Y & 0 \end{bmatrix} \ \begin{bmatrix} I & A \\ 0 & I \end{bmatrix} = \begin{bmatrix} 0 & 0 \\ Y & YA \end{bmatrix}. \]
Since $\begin{bmatrix} I & - A \\ 0 & I \end{bmatrix}$ is invertible with inverse $\begin{bmatrix} I &  A \\ 0 & I \end{bmatrix}$, it routinely follows that
\[
\begin{bmatrix} A Y - \alpha I & 0 \\ Y & -\alpha I\end{bmatrix}^k \text{\ \ \ \ \ and \ \ \ \ \ } \begin{bmatrix} -\alpha I & 0 \\ Y & Y A - \alpha I \end{bmatrix}^k \]
are similar.

An argument similar to that employed in \textsc{Step two} shows that
\[
\nul\, \begin{bmatrix} A Y - \alpha I & 0 \\ Y & -\alpha I\end{bmatrix}^k = \nul (AY - \alpha I)^k, \]
and that
\[
\nul\, \begin{bmatrix} -\alpha I & 0 \\ Y & Y A - \alpha I \end{bmatrix}^k = \nul \, (Y A - \alpha I)^k. \]
It follows that
\[
\nul\, (AY - \alpha I)^k = \nul\, (YA -\alpha I)^k, \]
and therefore
\begin{align*}
\nul \, (T + \alpha I)^k
	&= \nul(A Y + \alpha I)^k + \nul (Y A - \alpha I)^k \\
	&= \nul(Y A + \alpha I)^k + \nul (A Y - \alpha I)^k \\
	&= \nul(T - \alpha I)^k.
\end{align*}	
%
%
%%%%%%%%%%%%%%%%%%%%%%
%\bigskip
%%%%%%%%%%%%%%%%%%%%%%
%	
%$\bullet$\ \ \textsc{Step eight.}\ \ \ 	
%If $\dim\, \hilb < \infty$, we may write $T = Q \oplus T_1$, where $\sigma(Q) = \{ 0\}$ and $0 \not \in \sigma(T_1)$.    Then for $0 \ne \alpha$,
%\[
%\nul (T_1 - \alpha I)^k = \nul\, (T - \alpha I)^k = \nul\, (T+\alpha I)^k = \nul\, (T_1 + \alpha I)^k, \]
%whence $T_1 \sim - T_1$.    Since $Q$ is nilpotent (as it is quasinilpotent and $\dim\, \hilb < \infty$, it is clear by considering Jordan forms that $Q \sim - Q$.
%
%Thus $T = Q \oplus T_1 \sim -Q \oplus -T_1 = - T$, completing the proof.
%
\end{pf}

%%%%%%%%%%%%%%%%%%%%%%%%%%%%%%%%%%%%%%%%%%

\begin{cor} \label{cor2.09}
If $T \in \cniltwo$, then $\sigma_p(T) = \sigma_p(-T)$, and if $\dim\, \hilb < \infty$, then $T \sim -T$.
\end{cor}
	
%%%%%%%%%%%%%%%%%%%%%%%%%%%%%%%%%%%%%%%%%%

\bigskip

We require the following result of Barnes~\cite[Theorem~3]{Barnes1998}.

\begin{prop}\label{prop2.10} \emph{\textbf{(Barnes)}}\ \ \

Let $X$ and $Y$ be Banach spaces, and let $S\in \mathcal{B}(X,Y), R\in \mathcal{B}(Y,X)$, we have
\[\sigma_{ap}(RS)\setminus\{0\} =\sigma_{ap}(SR)\setminus \{0\}.\]
\end{prop}

%%%%%%%%%%%%%%%%%%%%%%%%%%%%%%%%%%%%%%%%%%

\begin{lem}\label{lem2.11}
Let $T \in \bofh$. Suppose that $T\in \mathfrak{c}(\ttt{nil}_2)$. Then $\sigma_{ap}(T)=\sigma_{ap}(-T)$.
\end{lem}

\begin{pf}
Let $T = [M, N]$, where $M^2 = 0 = N^2$, and write $M = \begin{bmatrix} 0 & A \\ 0 & 0 \end{bmatrix}$ and $N =\begin{bmatrix} W & X \\ Y & Z \end{bmatrix}$ relative to $\ker\, M \oplus (\ker\, M)^\perp$.  Then $T =\begin{bmatrix} AY &  AZ-WA\\ 0 & -YA \end{bmatrix}$.
We claim that
\[
\sigma_{ap}(T)\setminus\{0\}=(\sigma_{ap}(AY)\setminus\{0\})\cup (\sigma_{ap}(-YA)\setminus\{0\}). \]

$\bullet$\ \ \textsc{Step one.}\ \ \ 	Let $\lambda\in \sigma_{ap}(T)\setminus\{0\}$. Then there exists a sequence $\left( \begin{bmatrix} x_n\\ y_n \\ \end{bmatrix} \right)_{n=1}^{\infty}$of unit vectors
 such that
\[
\lim_n \begin{bmatrix} AY-\lambda &  AZ-WA\\ 0 & -YA-\lambda \end{bmatrix}\begin{bmatrix} x_n\\ y_n \\ \end{bmatrix} = \begin{bmatrix} 0\\ 0 \\ \end{bmatrix}.\]
If $\underset{n\in \mathbb{N}}{\inf}\|y_n\|>0$, it is easy to see that $\lambda\in \sigma_{ap}(-YA)$;
otherwise $\underset{n\in \mathbb{N}}{\inf}\|y_n\|=0$, and it is then easy to show that $\lambda\in \sigma_{ap}(AY)$.
Therefore,
\[\sigma_{ap}(T)\setminus\{0\}\subset(\sigma_{ap}(AY)\setminus\{0\})\cup (\sigma_{ap}(-YA)\setminus\{0\}).\]

$\bullet$\ \ \textsc{Step two.}\ \ \ 	
\begin{itemize}
	\item
	Let $0\neq \alpha\in \sigma_{ap}(AY)$, and $(x_n)_n \subset \ker M$ be a sequence of unit vectors such that $\lim_n \|(AY-\alpha)x_n\| = 0$. Then
    	\[
    	\lim_n \begin{bmatrix} AY-\alpha &  AZ-WA\\ 0 & -YA-\alpha \end{bmatrix}\begin{bmatrix} x_n\\ 0 \\ \end{bmatrix} = \begin{bmatrix} 0\\ 0 \\ \end{bmatrix},\]
    	which in turn implies that $\alpha\in \sigma_{ap}(T)$.
	\item
	Let $0\neq \beta\in \sigma_{ap}(-YA)$, and $\{y_n\}\subset (\ker M)^\perp$ be a sequence of unit vectors such that
	$\lim_n \|(-YA-\beta)y_n\| = 0$. Let $x_n=-\frac{1}{\beta}WAy_n$, $n\in \mathbb{N}$.
    	Define $z_n=\begin{bmatrix} x_n\\ y_n \\ \end{bmatrix}$, then $\|z_n\|\geq 1$, $n\in \mathbb{N}$.
    	We shall require the operators $L := MN$ and $R := NM$.
    	Then
    	\[
	\lim_n (R+\beta)z_n= \lim_n \begin{bmatrix} \beta &  WA\\ 0 & YA+\beta \end{bmatrix}\begin{bmatrix} x_n\\ y_n \\ \end{bmatrix} = \begin{bmatrix} 0\\ 0 \\ \end{bmatrix}.\]
	As $LR=RL=0$,
	\[
	\lim_n \beta Lz_n = \lim_n L[(R+\beta)z_n] = \begin{bmatrix} 0\\ 0 \\ \end{bmatrix}.\]
    	Since $\beta\neq 0$, it follows $\lim_n Lz_n = 0$.
    	Hence,
   	 \[
	 \lim_n (T-\beta)z_n = \lim_n (L-R-\beta)z_n = \lim_n Lz_n-(R+\beta)z_n = 0,\]
     	and thus $\beta\in \sigma_{ap}(T)$.
\end{itemize}
Therefore,
\[\sigma_{ap}(T)\setminus\{0\}\supset(\sigma_{ap}(AY)\setminus\{0\})\cup (\sigma_{ap}(-YA)\setminus\{0\}),\]
and the claim that $\sigma_{ap}(T)\setminus\{0\}=(\sigma_{ap}(AY)\setminus\{0\})\cup (\sigma_{ap}(-YA)\setminus\{0\})$ holds.

By Proposition~\ref{prop2.10},
\begin{align*}
\sigma_{ap}(T)\setminus\{0\}
	& =(\sigma_{ap}(AY)\setminus\{0\})\cup (\sigma_{ap}(-YA)\setminus\{0\}) \\
	&= (\sigma_{ap}(YA)\setminus\{0\})\cup (\sigma_{ap}(-AY)\setminus\{0\})\\
	&= \sigma_{ap}(-T)\setminus\{0\}.
\end{align*}	
This verifies that
$\sigma_{ap}(T)=\sigma_{ap}(-T)$.
\end{pf}

%%%%%%%%%%%%%%%%%%%%%%%%%%%%%%%%%%%%%%%%%%

\bigskip

Let $\hilb$ be an infinite-dimensional, separable Hilbert space and $\cM_1$, $\cM_2$ be infinite-dimensional subspaces of $\hilb$ such that $\cM_2 = \cM_1^\perp$.  Given $T \in \bofh$, write $T = \begin{bmatrix} T_{11} & T_{12} \\ T_{21} & T_{22} \end{bmatrix}$ relative to the decomposition $\hilb = \cM_1 \oplus \cM_2$.   Since $\cM_1, \cM_2$ and $\hilb$ are isomorphic as Hilbert spaces, this allows us to identify $\bofh = \cB(\cM_1 \oplus \cM_2)$ with $\cB(\hilb \oplus \hilb)$.  In other words, up to Hilbert space isomorphism and the usual identification of $\cB(\hilb \oplus \hilb)$ with $\bbM_2(\bofh)$, there is no harm in assuming that $T_{ij} \in \bofh$, $1 \le i, j \le 2$.   We shall employ this standard device in the proof of item (b) of the next proposition.

\bigskip

\begin{prop}\label{prop2.12}
Let $T \in \bofh$. Suppose that $T\in \mathfrak{c}(\ttt{nil}_2)$. Then
\begin{enumerate}
	\item[(a)]
	every connected component of $\sigma(T)$ intersects $\sigma(-T)$;
	\item[(b)]
	every connected component of $\sigma_e(T)$ intersects $\sigma_e(-T)$.
\end{enumerate}
\end{prop}

\begin{pf}
Let $T = [M, N]$ where $M^2 = 0 = N^2$. Writing $M = \begin{bmatrix} 0 & A \\ 0 & 0 \end{bmatrix}$ and $N =\begin{bmatrix} W & X \\ Y & Z \end{bmatrix}$ relative to $\ker\, M \oplus (\ker\, M)^\perp$,
and so $T =\begin{bmatrix} AY &  AZ-WA\\ 0 & -YA \end{bmatrix}$.

\begin{enumerate}
	\item[(a)]
	Let $\Omega$ be a component of $\sigma(T)$. Then $\Omega\cap \partial \sigma(T)\neq \varnothing$.
Since $\partial \sigma(T)\subset \sigma_{ap}(T)$(\cite{Conway1990} VII Proposition 6.7), it follows that $\Omega\cap\sigma_{ap}(T)\neq \varnothing$.
By Lemma \ref{lem2.11}, $\Omega\cap\sigma_{ap}(-T)\neq \varnothing$. In particular,
$\Omega\cap\sigma(-T)\neq \varnothing$.
	\item[(b)]
	The case where $\dim\, \hilb < \infty$ is handled by Theorem~\ref{thm2.08}.
	
	When $\dim\, \hilb = \infty$, the fact that $M^2 = 0$ implies that $\dim\, \ker\, M = \infty$.    If $\dim\, (\ker\, M)^\perp < \infty$, then both $M$ and (therefore) $T$ have finite rank, whence $\sigma_e(T)=\sigma_e(-T)=\{0\}$.   As such, we may assume without loss of generality that $\dim \ker\, M=\dim (\ker\, M)^\perp=\dim\, \hilb = \infty$.
	
	As per the comments preceding this Proposition, we think of $A, W, Y$ and $Z$ as elements of $\bofh$.

Let $\pi$ be the quotient map $\pi$ from $\mathcal{B}(\cH)$ onto the Calkin algebra $\mathcal{B}(\cH)/\mathcal{K}(\cH)$.
 Then
 \[t:=\pi\otimes \textup{id}_{\mathbb{M}_2}(T)=\begin{bmatrix} \pi(A)\pi(Y) &  \pi(AZ-WA)\\ 0 & -\pi(Y)\pi(A) \end{bmatrix}.\]
 Let $\rho$ be a unital faithful $C^*$-representation of $C^*(\pi(A),\pi(Y),\pi(Z),\pi(W))$ to $\mathcal{B}(\cL)$,
 where $\cL$ is a separable Hilbert space.
 Then
 \[\tilde{T}:=(\rho\otimes \textup{id}_{\mathbb{M}_2})(t)=\begin{bmatrix} (\rho\circ\pi)(A)\cdot(\rho\circ\pi)(Y) &  (\rho\circ\pi)(AZ-WA)\\ 0 & -(\rho\circ\pi)(Y)\cdot(\rho\circ\pi)(A) \end{bmatrix}.\]
Moreover, $\sigma_e(T)=\sigma(t)=\sigma(\tilde{T})$ and $\tilde{T}\in \mathfrak{c}(\ttt{nil}_2)$.
The proof is completed by applying (a).
\end{enumerate}
\end{pf}

%%%%%%%%%%%%%%%%%%%%%%%%%%%%%%%%%%%%%%%%%%

%%%%%%%%

%%%%%%%%%%%%%%%%%%%%%%%%%%%%%%%%%%%%%%%%%%
%%%%%%%%%%%%%%%%%%%%%%%%%%%%%%%%%%%%%%%%%%
% SECTION THREE
%%%%%%%%%%%%%%%%%%%%%%%%%%%%%%%%%%%%%%%%%%
%%%%%%%%%%%%%%%%%%%%%%%%%%%%%%%%%%%%%%%%%%

%%%%%%%%

%%%%%%%%%%%%%%%%%%%%%%%%%%%%%%%%%%%%%%%%%%

\vskip 1.5 cm

\section{$\cniltwo$ in the finite-dimensional setting} \label{section3}

%%%%%%%%

%%%%%%%%%%%%%%%%%%%%%%%%%%%%%%%%%%%%%%%%%%

\subsection{} \label{sec3.01}
In this section we shall characterise the set $\cniltwo$ in the case where $\dim\, \hilb = n < \infty$.

Suppose that $T \in \cB(\bbC^n)$.  By considering the Jordan form of $T$, we find that we may write $T \sim B \oplus J$, where $B$ is invertible and $J$ is nilpotent.

Proposition~\ref{prop2.02} above asserts that $T \in \cniltwo$ if and only if $B \in \cniltwo$ and $J \in \cniltwo$.

%%%%%%%%%%%%%%%%%%%%%%%%%%%%%%%%%%%%%%%%%%

\begin{cor} \label{cor3.02}
Let $n \in \bbN$, and suppose that $T \in \cB(\bbC^n)$ lies in $\cniltwo$.    Then $T \sim -T$.
\end{cor}

\smallskip

\begin{pf}
Let us write the Jordan form of $T$ as $B \oplus J$, where $B$ is invertible and $J$ is nilpotent.   By  Proposition~\ref{prop2.02}, $B \in \cniltwo$ and $J \in \cniltwo$.    Now, by Proposition~\ref{prop2.03} and the comment which follows it, $B \sim -B$.    Since $J$ is nilpotent and acts on a finite-dimensional space, $J \sim - J$, from which it easily follows that $T \sim -T$.
\end{pf}

%%%%%%%%%%%%%%%%%%%%%%%%%%%%%%%%%%%%%%%%%%

The converse to Corollary~\ref{cor3.02} holds for normal operators (Theorem~\ref{thm2.06}).   In general, however, it is false.   To see this, we first appeal to a result from ~\cite{RadjaviRosenthal2002}  which we shall also require in the sequel.

%%%%%%%%%%%%%%%%%%%%%%%%%%%%%%%%%%%%%%%%%%

\begin{thm} \label{thmRR2002}  \cite[Theorem 1]{RadjaviRosenthal2002}.
If $A$ and $B$ are quadratic operators on a Banach space and $A B - B A$ is nilpotent, then $\{ A, B \}$ is triangularizable.
\end{thm}

%%%%%%%%%%%%%%%%%%%%%%%%%%%%%%%%%%%%%%%%%%

\bigskip

Our characterisation of $\cniltwo$ in Theorem~\ref{thm3.15} will subsume the next result. Having said that, the result is elementary and can easily be used to exhibit a large number of nilpotent operators which do not lie in $\cniltwo$.  In particular, it may be used to provide a large number of examples of operators $T \sim -T$ which do not lie in $\cniltwo$.

\begin{cor} \label{cor3.04}
Let $J \in \cB(\bbC^n)$ is a nilpotent operator in $\cniltwo$, then $J^{\lfloor \frac{n+1}{2} \rfloor} = 0$.   In particular, if $4 \le n$ and $J_n \in \cB(\bbC^n)$ is the usual nilpotent Jordan nilpotent of order $n$, then $J_n \sim -J_n \not \in \cniltwo$.
\end{cor}

\begin{pf}
Let $J \in \cB(\bbC^n)$ be nilpotent, and suppose that $J \in \cniltwo$, say $J = [M, N]$ where $M^2 = 0 = N^2$.   Let $L := MN$ and $R := NM$.   Clearly $L R = 0 = RL$, from which we deduce that $J^k = (L - R)^k = L^k + (-1)^k R^k$ for all $k \ge 1$.

Since $J = [M, N]$ is nilpotent, by Theorem~\ref{thmRR2002}, there exists an orthonormal basis $\{ e_k\}_{k=1}^n$ for $\bbC^n$ relative to which the matrices $[M]$ and $[N]$ corresponding to $M$ and $N$ respectively are upper-triangular.   In fact, since $M^2 = 0 = N^2$, $[M]$ and $[N]$ are strictly upper-triangular.
But the algebra $\fT_n^\circ \subseteq \bbM_n(\bbC)$ of strictly upper-triangular $n \times n$ complex matrices is a nil algebra and clearly $A_1 A_2 \cdots A_n = 0$ for all choices of $A_k \in \fT_n^\circ$, $1 \le k \le n$.   It follows that $L^k = 0 = R^k$, and thus $J^k = 0$, if $k \ge \lfloor \frac{n+1}{2} \rfloor$.

The second sentence in the statement of the Corollary is an immediate consequence of the first.
\end{pf}

%%%%%%%%%%%%%%%%%%%%%%%%%%%%%%%%%%%%%%%%%%

\subsection{} \label{sec3.05}
As previously observed, an element $T \in \cB(\bbC^n)$ lies in $\cniltwo$ if and only if $B \in \cniltwo$ and $J \in \cniltwo$, where $T \sim B \oplus J$, where $B$ is invertible and $J$ is nilpotent.  Proposition~\ref{prop2.03} tells us that $B \in \cniltwo$ if and only if $B \sim A_0 \oplus - A_0$ for some invertible operator $A_0$.
We now turn our attention to the problem of determining which nilpotent operators in $\cB(\bbC^n)$ lie in $\cniltwo$.   To do so, we shall need to calculate the form of commutators and anti-commutators of specific Jordan forms.   Similar results have appeared in the literature previously, and we draw the reader's attention to~\cite{OzturkKaptanoglu2010, Oblak2012, Khatami2013}.

\bigskip

We begin with a surprisingly useful result.  But before that, we remind the reader that if $\hilb$ is a Hilbert space, $T \in \bofh$ and $n \in \bbN$, then $\hilb^{(n)} := \hilb \oplus \hilb \oplus \cdots \oplus \hilb$ and $T^{(n)} := T \oplus T \oplus \cdots \oplus T$, where in each case we have $n$ summands.

%%%%%%%%%%%%%%%%%%%%%%%%%%%%%%%%%%%%%%%%%%

\begin{prop} \label{prop3.06}
Let $\mu \in \bbN$.    Then $J_2^{(\mu)} \in \cniltwo$ if and only if $\mu$ is even.
\end{prop}

\begin{pf}
First note that $J_2^{(2)} = J_2 \oplus J_2 \sim J_2 \oplus - J_2 \simeq [M, N]$, where
\[
M = \begin{bmatrix} 0 & 1 & 0 & 0 \\ 0 & 0 & 0 & 0 \\ 0 & 0 & 0 & 1 \\ 0 & 0 & 0 & 0 \end{bmatrix}, \ \ \ \ \
N = \begin{bmatrix} 0 & 0 & 0 & 0 \\ 0 & 0 & 1 & 0 \\ 0 & 0 & 0 & 0 \\ 0 & 0 & 0 & 0 \end{bmatrix}. \]
Thus if $\mu = 2 \kappa$ for some $\kappa \ge 1$, then
\[
J_2^{(\mu)} = J_2^{(2 \kappa)} \sim (J_2 \oplus - J_2)^{(\kappa)} \simeq [M, N]^{(\kappa)} = [M^{(\kappa)}, N^{(\kappa)}]. \]
Since $M^{(\kappa)}$ and $N^{(\kappa)}$ are nilpotent of order two, we see that $J_2^{(\mu)} \in \cniltwo$ when $\mu$ is even.

\bigskip

Next we will show that if $0 \le \kappa$ is an integer, then $J_2^{(2\kappa+1)}$ is not in $\cniltwo$. Otherwise, there exist square-zero operators $M,N \in \cB(\bbC^{2(2\kappa+1)})$ such that
\[
J := \begin{bmatrix}0& I_{2\kappa +1} \\ 0& 0 \end{bmatrix} = [M, N].\]
As $M J = - J M$ and $N J = - J N$, relative to the same decomposition of $\hilb$, $M$ and $N$ admit the following decompositions:
\[M=\begin{bmatrix}M_0& M_1 \\ 0 &  -M_0  \end{bmatrix},~~N=\begin{bmatrix}N_0& N_1 \\ 0 & -N_0 \end{bmatrix},\]
where $M_i,N_i\in \cB(\bbC^{2\kappa+1})$, $i=0,1$.

Since $M^2 = 0$, it follows that $M_0^2 = 0$, and so we can find an invertible operator $S \in \cB(\bbC^{2 \kappa + 1})$ and $0 \le s \le \kappa$ such that $S^{-1} M_0 S = \begin{bmatrix} 0 & I_s \\ 0 & 0 \end{bmatrix} \oplus 0^{(2 \kappa + 1 - 2 s)}$.  Observe that
\[
(S\oplus S)^{-1} J (S \oplus S) = J, \]
and thus by replacing $M$ and $N$ by $(S\oplus S)^{-1} M (S \oplus S)$ and $(S \oplus S)^{-1} N (S \oplus S)$, we may assume without loss of generality that (relative to a decomposition of $\bbC^{4 \kappa +2} = \hilb_1 \oplus \hilb_2 \oplus \hilb_3$)
\[
M_0 = \begin{bmatrix} 0 & I_s & 0 \\ 0 & 0 & 0 \\ 0 & 0 & 0 \end{bmatrix} = \begin{bmatrix} 0 & I_s \\ 0 & 0 \end{bmatrix} \oplus 0^{(2 \kappa + 1 - 2 s)}. \]
Obviously, $J \not = 0$ implies that $s \ge 1$.   Furthermore,  $N^2 = 0$  implies that $N_0^2 = 0$, and the equations $M^2 = 0$ and $J = [M, N]$ imply that
\begin{itemize}
	\item{}
	$M_0 M_1=M_1 M_0$;
	\item{}
	$M_0 N_0 = N_0 M_0$; and
	\item{}
	$M_0N_1-M_1N_0-N_0M_1+N_1M_0=I_{2\kappa +1}.$
\end{itemize}

The first two equations  imply that relative to the decomposition $\bbC^{4 \kappa + 2} = \hilb_1 \oplus \hilb_2 \oplus \hilb_3$ from above,  $M_1$ and $N_0$  have the following forms:
\[M_1=\begin{bmatrix} A_0& A_1& B\\ 0& A_0& 0\\ 0& C& D\\ \end{bmatrix},~~N_0=\begin{bmatrix} E_0& E_1 & F\\ 0& E_0& 0\\0 & G& K\\ \end{bmatrix},
\]
where $A_i,E_i\in \mathbb{M}_{s}(\bbC)$, $i=0,1$, and $D, K\in \cB(\bbC^{2\kappa+1 - 2s})$.  As $N_0$ is square-zero, so is $K$.

A straightforward calculation based upon the third equation above shows that the compression of $I_{2 \kappa + 1}$ to $\hilb_3$ is nothing more than
\[
I_{2 \kappa + 1  - 2s} = - (D K + K D). \]
Now $K \in \cB(\bbC^{2\kappa+1 - 2s })$ is a square-zero operator, and so $K$ is similar to $J_2^{(t)} \oplus 0^{(2 \kappa + 1 - 2s -2t)}$.    Crucially, we note that $2\kappa+1 - 2s$ is an odd number, and therefore $(2 \kappa + 1) - 2 s - 2 t \ge 1$.   Noting that $I_{2 \kappa + 1 - 2s}$ is invariant under conjugation by invertible elements,   we may choose a basis for $\bbC^{2 \kappa + 1 - 2 s}$ with respect to which the last column and last row of the corresponding matrix for $K$ are both zero.  It follows that the $(2 \kappa+ 1 - 2s, 2 \kappa + 1 - 2s)$ entry of  $I_{2 \kappa + 1 - 2s} = -(D K + K D)$ is zero, which is clearly false.

This contradiction implies that $J_2^{(2\kappa+1)} \simeq J$ is not in $\cniltwo$.
\end{pf}

%%%%%%%%%%%%%%%%%%%%%%%%%%%%%%%%%%%%%%%%%%

Our next goal is to establish that a nilpotent element $J$ of $\cB(\bbC^n)$ lies in $\cniltwo$ if and only if there exists a nilpotent operator $Q \in \cB(\bbC^n)$ such that $J = Q^2$.    To do this, we first require the a characterisation of those nilpotent matrices admitting a square root.   Fortunately, we may appeal to a result of Borwein and Richmond which establishes the conditions under which an arbitrary $n \times n$ complex matrix admits a square root.   We apply their result to the case of nilpotent matrices.

\begin{thm}  \emph{(\textbf{Borwein and Richmond})\cite[Theorem 1]{BorweinRichmond1984}} \label{thm3.07}

Let $2 \le n \in \bbN$ and let $J \in \cB(\bbC^n)$ be a nilpotent operator.   Choose integers $1 \le \gamma$ and $\alpha_1 \ge \alpha_2 \ge \cdots \ge \alpha_\gamma$  such that
\[
J \sim J_{\alpha_1} \oplus J_{\alpha_2} \oplus \cdots J_{\alpha_\gamma}. \]
If $\gamma$ is odd, we define $\alpha_{\gamma+1} = 0$.
Then $J$ is the square of a (necessarily nilpotent) operator $Q$ if and only if
\[
\alpha_{2 j - 1} - \alpha_{2 j} \le 1, \ \ \ 1 \le j \le \lfloor \frac{\gamma + 1}{2} \rfloor. \]
\end{thm}

%%%%%%%%%%%%%%%%%%%%%%%%%%%%%%%%%%%%%%%%%%

The precise reformulation of their result which we shall use below is as follows.

\begin{thm} \label{thm3.08}

Let $2 \le n \in \bbN$ and suppose that $J \in \cB(\bbC^n)$ is nilpotent.   We choose integers $1 \le d$, $\mu_1, \mu_2, \ldots, \mu_d$ and $\kappa_1 \ge \kappa_2 \ge \cdots \ge \kappa_d \ge 1$ such that
\[
J \sim J_{\kappa_1}^{(\mu_1)} \oplus J_{\kappa_2}^{(\mu_2)} \oplus \cdots \oplus J_{\kappa_d}^{(\mu_d)}. \]
The following conditions are equivalent.
\begin{enumerate}
	\item[(a)]
	$J = Q^2$ for some nilpotent operator $Q \in \cB(\bbC^n)$.
	\item[(b)]
	\begin{itemize}	
		\item{}
		If $1 \le b \le d-1$ and $\kappa_b - \kappa_{b+1} \ge 2$, then $\sum_{i=1}^b \mu_i$ is an even integer; and
		\item{}
		if $b = d$ and $\kappa_b \ge 2$, then $\sum_{i=1}^b \mu_i$ is an even integer.
	\end{itemize}
\end{enumerate}	
\end{thm}

%%%%%%%%%%%%%%%%%%%%%%%%%%%%%%%%%%%%%%%%%%

\subsection{} \label{sec3.09}
The keys to our characterisation of $\cniltwo$ in the finite-dimensional setting are Theorem~\ref{thm3.08} above, together with our previous observation that if $J = [M, N]$ where $M^2 = 0 = N^2$, then $M J = - J M $ and $N J = - J N$.

%%%%%%%%%%%%%%%%%%%%%%%%%%%%%%%%%%%%%%%%%%

In light of these two anti-commutation relations, the following Lemma will prove useful to us below.  The proof is an elementary computation.  The reader may wish to compare this result as well as our results concerning the structure of $M$ and $N$ anti-commuting with $J$ as above with the results of Oblak~\cite{Oblak2012} on the structure of nilpotent operators commuting with a fixed nilpotent matrix.

\begin{lem} \label{lem3.10}
Let $m, n \ge 1$ be integers and suppose that $D = [d_{ij}] \in \bbM_{m\times n} (\bbC)$.   Suppose furthermore that $J_m D = - D J_n$, and set $\alpha_j = d_{1 j}$, $1 \le j \le n$.
\begin{enumerate}
	\item[(a)]
	If $2 \le m = n$, then
	\[
	D =
	\begin{bmatrix}
		\alpha_1  & \alpha_2 & \alpha_3 & \cdots  & \alpha_{n-2} & \alpha_{n-1} & \alpha_{n}  \\
		0 & -\alpha_1 & -\alpha_2 &  \cdots  &  \cdots & -\alpha_{n-2} & - \alpha_{n-1}  \\
		\vdots & & & & & &  \vdots \\
		0 & 0 & 0 & \cdots &  \cdots & (-1)^{n-2} \alpha_1  & (-1)^{n-2}\alpha_2  \\
		0 & 0 & 0 & \cdots & \cdots & 0 & (-1)^{n-1} \alpha_1 \\
	\end{bmatrix}. \]
	\item[(b)]
	If $1 \le m < n$, then
	\[
	D =
	\begin{bmatrix}
		0 & \cdots & 0 &  \alpha_{n-m + 1} & \alpha_{n-m+2} & \cdots & \cdots &  \alpha_n \\
		0 & \cdots & 0 &  0 & - \alpha_{n-m+1} & - \alpha_{n-m+2} & \cdots & -\alpha_{n-1} \\
		\vdots & & & & & & & \vdots \\
		 0 & \cdots & 0 &  \cdots &   & \cdots &  (-1)^{m-2} \alpha_{n-m+1} & (-1)^{m-2}  \alpha_{n-m+2} \\
		 0 & \cdots & 0 &  0 & \cdots &  \cdots  & 0 & (-1)^{m-1} \alpha_{n-m+1} \\
	\end{bmatrix}. \]
	In particular, the first $n-m$ columns of $D$ are zero, and the last $m$ columns of $D$ form an upper-triangular matrix.
	\item[(c)]
	If $1 \le n < m$, then
	\[
	D=
	\begin{bmatrix}
		\alpha_1  & \alpha_2 & \alpha_3 & \cdots  & \alpha_{n-2} & \alpha_{n-1} & \alpha_{n}  \\
		0 & -\alpha_1 & -\alpha_2 &  \cdots  &  \cdots & -\alpha_{n-2} & - \alpha_{n-1}  \\
		\vdots & & & & & &  \vdots \\
		0 & 0 & 0 & \cdots &  \cdots & (-1)^{n-2} \alpha_1  & (-1)^{n-2}\alpha_2  \\
		0 & 0 & 0 & \cdots & \cdots & 0 & (-1)^{n-1} \alpha_1 \\
		0 & 0 & 0 & \cdots & 0 & 0 & 0\\
		\vdots &  &  & \cdots &  &  & \vdots \\
		0 & 0 & 0 & \cdots & 0 & 0 & 0\\	
	\end{bmatrix}. \]
	In particular, the last $m-n$ rows of $D$ are zero, and the first $n$ rows of $D$ form an upper-triangular matrix.
\end{enumerate}
\end{lem}

%%%%%%%%%%%%%%%%%%%%%%%%%%%%%%%%%%%%%%%%%%

\begin{rem} \label{rem3.11}
Although the forms described in detail in the above Lemma are in themselves interesting, for our purposes, we shall only require the following observation.  For reasons which will become apparent below, we now replace $n$ and $m$ above by $\kappa_{i_1}$ and $\kappa_{i_2}$ respectively, and we denote the orthogonal bases for $\bbC^{\kappa_{i_1}}$ and $\bbC^{\kappa_{i_2}}$ by $\{ e_{(i_1, j_1, s)} : 1 \le s \le \kappa_{i_1} \}$ and $\{ e_{(i_2, j_2, t)} : 1 \le t \le \kappa_{i_2}\}$.  We have just shown that if $\max(\kappa_{i_1}, \kappa_{i_2}) \ge 2$, and if $D = [d_{ij}] \in \bbM_{\kappa_{i_2} \times \kappa_{i_1}}(\bbC)$ satisfies $J_{\kappa_{i_2}} D = - D J_{\kappa_{i_1}}$, then:
\begin{enumerate}
	\item[(\textsc{i})]
	$\langle D e_{(i_1, j_1, t)}, e_{(i_2, j_2, \kappa_{i_2})} \rangle \ne 0$ implies that $\kappa_{i_2} \le \kappa_{i_1}$ and $t= \kappa_{i_1}$.
	\item[(\textsc{ii})]
	$\langle D e_{(i_1, j_1, t)}, e_{(i_2, j_2, \kappa_{i_2} -1)} \rangle \ne 0$ implies that $\kappa_{i_2}  - 1 \le \kappa_{i_1}$ and either
	\begin{itemize}
		\item{}
		$\kappa_{i_2} - 1 = \kappa_{i_1}$, and $t = \kappa_{i_1}$, or
		\item{}
		$\kappa_{i_2} \le \kappa_{i_1}$, and $t \in \{ \kappa_{i_1} -1, \kappa_{i_1}\}$.
	\end{itemize}	
\end{enumerate}	

\end{rem}

%%%%%%%%%%%%%%%%%%%%%%%%%%%%%%%%%%%%%%%%%%

The following result came as quite a surprise to the authors.

\begin{thm} \label{thm3.12}
Let $2 \le n \in \bbN$ and $J \in \cB(\bbC^n)$ be a nilpotent operator.   If $J \in \cniltwo$, then $J$ admits a square root.
\end{thm}	

\begin{pf}
\noindent{\textsc{Step 1.}}

The uniqueness of the Jordan canonical form of $J$ (up to permutation of the irreducible summands) allows us to find integers $1 \le d, \mu_1, \mu_2, \ldots, \mu_d$ and $\kappa_1 \ge \kappa_2 \ge \cdots \ge \kappa_d \ge 1$ such that
\[
J := J_{\kappa_1}^{(\mu_1)} \oplus J_{\kappa_2}^{(\mu_2)} \oplus \cdots \oplus J_{\kappa_d}^{(\mu_d)} \]
relative to the decomposition $\bbC^n = \hilb_1^{(\mu_1)} \oplus \hilb_2^{(\mu_2)} \oplus \cdots \oplus \hilb_d^{(\mu_d)}$.  Let $\hilb_{(i,j)}$ denote the $j^{th}$ copy of $\hilb_i \simeq \bbC^{\kappa_i}$ in $\hilb_i^{(\mu_i)}$, $1 \le j \le \mu_i$, and let $\{ e_{(i, j, s)}: 1 \le s \le \kappa_i\}$ denote the orthonormal basis of $\hilb_{(i,j)}$ that gives rise to the standard super-diagonal form for the $j^{th}$ copy of $J_{\kappa_i}$ in $J_{\kappa_i}^{(\mu_i)}$.   Also, let $P_{(i,j)}$ denote the orthogonal projection of $\bbC^n$ onto $\hilb_{(i,j)}$,  $1 \le j \le \mu_i$, $1 \le i \le d$.

We have that $J = [M, N]$, that $M^2 = 0 = N^2$ and thus $J M = - M J$, $J N  = - N J$.   Observe that
\[
P_{(i_2, j_2)} M P_{(i_1,j_1)} \in \cB(\hilb_{(i_1, j_1)}, \hilb_{(i_2,j_2)}) \simeq \bbM_{\kappa_{i_2} \times \kappa_{i_1}}(\bbC). \]
A routine calculation ensuing from the fact that $J M = -  M J$ yields
\[
J_{\kappa_{i_2}} (P_{(i_2, j_2)} M P_{(i_1, j_1)}) = - (P_{(i_2, j_2)} M P_{(i_1,j_1)}) J_{\kappa_{i_1}}. \]

%(Perhaps unsurprisingly, the operators $(P_{(i,j)} M P_{(s, t)}$ will play the role of $D$ in Lemma~\ref{lem2.12} and Remark~\ref{rem2.13} above.)

%%%%%%%%%%%%%%%%%%%%%

\bigskip
\noindent{\textsc{Step 2.}}

Suppose that $1 \le b < d$ and that $\kappa_b - \kappa_{b+1} \ge 2$.  We define four spaces $\fM_\ell$, $1 \le \ell \le 4$.

\begin{itemize}
	\item{}
	We define $\fM_3 := \mathrm{span} \{ e_{(i, j, \kappa_i)} : 1 \le i \le b, 1 \le j \le \mu_i \}$.  In other words, we are choosing the last basis vector $e_{(i,j,\kappa_i)}$ from each copy $\hilb_{(i,j)}$ of the space upon which each $J_{\kappa_i}$ acts, \emph{up to and including all copies of $\hilb_{(b, j)}$, $1 \le j \le \mu_b$}, but no copies occurring after that.
	\item{}
	We define $\fM_2 := \mathrm{span} \{ e_{(i, j, \kappa_i - 1)} : 1 \le i \le b, 1 \le j \le \mu_i \}$.  In other words, we are choosing the \emph{second to last} basis vector $e_{(i,j,\kappa_i -1)}$from each copy $\hilb_{(i,j)}$ of the space upon which each $J_{\kappa_i}$ acts, \emph{up to and including all copies of $\hilb_{(b, j)}$, $1 \le j \le \mu_b$}, but no copies occurring after that.
	\item{}
	We define $\fM_1 := \mathrm{span} \{ e_{(i,j, s)} : 1 \le i \le b, 1 \le j \le \mu_i, 1 \le s \le \kappa_i -2 \}$.     This is the space
	\[
	\left( \hilb_1^{(\mu_1)} \oplus \cdots \oplus \hilb_b^{(\mu_b)} \right) \ominus (\fM_2 \oplus \fM_3). \]
	\item{}
	We define
	\[
	\fM_4 := \bbC^n \ominus (\fM_1 \oplus \fM_2 \oplus \fM_3) =  \hilb_{b+1}^{(\mu_{b+1})} \oplus \hilb_{b+2}^{(\mu_{b+2})} \oplus \cdots \oplus \hilb_d^{(\mu_d)}. \]
\end{itemize}	

%%%%%%%%%%%%%%%%%%%%%

\bigskip
\noindent{\textsc{Step 3.}}

Write $M = [M_{i,j}]$ relative to the decomposition
\[
\bbC^n = \fM_1 \oplus \fM_2 \oplus \fM_3 \oplus \fM_4, \]
and note that with respect to this decomposition, $J$ is of the form
\[
J =
\begin{bmatrix}
	Q_{1,1} & Q_{1,2} & 0 & 0 \\
	0 & 0 & I_\mu & 0 \\
	0 & 0 & 0 & 0 \\
	0 & 0 & 0 & Q_{4,4}
\end{bmatrix}, \]
where $\mu := \sum_{i=1}^b  \mu_i$.
Our immediate goal is to show that the equation $J M = - M J$ implies that $M_{2,1} = M_{2, 4} = M_{3, 1} = M_{3, 2} = M_{3, 4} = 0$, so that $M$ is of the form
\[
\begin{bmatrix}
	M_{1, 1} & M_{1, 2} & M_{1, 3} & M_{1,4} \\
	0 & M_{2, 2} & M_{2,3} & 0 \\
	0 & 0 &  M_{3, 3} & 0 \\
	M_{4, 1} & M_{4, 2} & M_{4, 3} & M_{4, 4}
\end{bmatrix}. \]
(In fact, we can do better than this.  For example, it can be shown that $M_{3, 3} = - M_{2, 2}$, and we can better specify other matrix entries of $M$, but for our purposes the stated form of $M$ will suffice.)

To that end, choose an arbitrary basis vector $e_{(i_1, j_1, s_1)}$ for $\fM_1$;  i.e. $1 \le i_1 \le b$, $1 \le j_1 \le \mu_{i_1}$ and $1 \le s_1 \le \kappa_{i_1} - 2$.   Next, fix an arbitrary basis vector $e_{(i_2, j_2, \kappa_{i_2} - 1)}$ for $\fM_2$, so that $1 \le i_2 \le b$, $1 \le j_2 \le \mu_{i_2}$.   Note that $e_{(i_2, j_2, \kappa_{i_2})}$ is  an arbitrary basis vector for $\fM_3$.

Then
\[
\langle M e_{(i_1, j_1, s_1)}, e_{(i_2, j_2, \kappa_{i_2}-1)} \rangle = \langle P_{(i_2, j_2)} M P_{(i_1, j_1)} e_{(i_1, j_1, s_1)}, e_{(i_2, j_2, \kappa_{i_2}-1)} \rangle, \]
and
\[
\langle M e_{(i_1, j_1, s_1)}, e_{(i_2, j_2, \kappa_{i_2})} \rangle = \langle P_{(i_2, j_2)} M P_{(i_1, j_1)} e_{(i_1, j_1, s_1)}, e_{(i_2, j_2, \kappa_{i_2})} \rangle. \]

But $P_{(i_2, j_2)} M P_{(i_1, j_1)} \in \cB(\hilb_{(i_1, j_1)}, \hilb_{(i_2, j_2)}) \simeq \bbM_{\kappa_{i_2} \times \kappa_{i_1}}(\bbC)$, and so we may think of\linebreak $P_{(i_2, j_2)} M P_{(i_1, j_1)}$ as playing the role of $D$ in Lemma~\ref{lem3.10} and Remark~\ref{rem3.11} above. By that Remark,

\begin{itemize}
	\item{}
	if $\langle (P_{(i_2, j_2)} M P_{(i_1, j_1)}) e_{(i_1, j_1, s_1)}, e_{(i_2, j_2, \kappa_{i_2})}\rangle \ne 0$, then $\kappa_{i_2} \le \kappa_{i_1}$ and $s_1 = \kappa_{i_1}$. The condition that $\kappa_{i_2} \le \kappa_{i_1}$ implies that $i_1 \le i_2 \le b$, and so we may conclude that $M_{3, 4} = 0$.    The condition that $s_1 = \kappa_{i_1}$ then implies that  $e_{(i_1, j_1, s_1)} = e_{(i_1, j_1, \kappa_{i_1})}  \in \fM_3$.  In other words, $M_{3, 1} = M_{3, 2}  = 0$.
	\item{}
	Next, suppose that $\langle P_{(i_2, j_2)} M P_{(i_1, j_1)} e_{(i_1, j_1, s_1)}, e_{(i_2, j_2, \kappa_{i_2}-1)} \rangle \ne 0$.   By Remark~\ref{rem3.11}(\textsc{ii}), either $\kappa_{i_2} - 1 \le \kappa_{i_1}$.     The hypotheses that $i_2 \le b$ and  $\kappa_b - \kappa_{b+1} \ge 2$ then imply that $i_1 \le i_2$, whence $M_{2,4} = 0$.  The same remark also shows that $s_1 \in \{ \kappa_{i_1} - 1, \kappa_{i_1}\}$, implying that $e_{(i_1, j_1, s_1)} \in \fM_2 \oplus \fM_3$.   In other words, $M_{2,1} = 0$.
\end{itemize}
	
Hence $M$ has the desired form.  Since $M^2 = 0$, we conclude that $\begin{bmatrix} M_{2, 2} & M_{2, 3} \\ 0 & M_{3, 3} \end{bmatrix}$ is nilpotent of order two.

Since $J N = - N J$, the same argument implies that
\[
N =
\begin{bmatrix}
	N_{1,1} & N_{1, 2} & N_{1, 3} & N_{1, 4} \\
	0 & N_{2, 2} & N_{2, 3} & 0 \\
	0 & 0 & N_{3, 3} & 0 \\
	N_{4, 1} & N_{4, 2} & N_{4, 3} & N_{4, 4}
\end{bmatrix}, \]
and as before, since $N^2 = 0$, we conclude that $\begin{bmatrix} N_{2, 2} & N_{2, 3} \\ 0 &  N_{3, 3} \end{bmatrix}$ is nilpotent of order two.

%%%%%%%%%%%%%%%%%%%%%

\bigskip
\noindent{\textsc{Step 4.}}

Now $J = [M, N]$.   By comparing the compression of $J$ and of $[M, N]$ to $\fM_2 \oplus \fM_3$, we find that
\[
\begin{bmatrix} 0 & I_\mu \\ 0 & 0 \end{bmatrix} = \left[ \begin{bmatrix} M_{2, 2} & M_{2, 3} \\ 0 &  M_{3, 3} \end{bmatrix}, \begin{bmatrix} N_{2, 2} & N_{2, 3} \\ 0 &  N_{3, 3} \end{bmatrix} \right] \in \cniltwo.\]
By Proposition~\ref{prop3.06}, $\mu$ is even.  On the other hand,
\[
\mu = \dim\, \fM_3 = \dim\, \fM_2 = \sum_{i = 1}^b \mu_i. \]
In other words, $\mu$ counts exactly how many operators there are up to and including the $\mu_b^{th}$ copy of $J_{\kappa_b}$.

%%%%%%%%%%%%%%%%%%%%%

\bigskip
\noindent{\textsc{Step 5.}}

There remains to prove that if $\kappa_d \ge 2$, then $\sum_{i=1}^d \mu_i$ is even.  The proof is a routine adaptation of the previous case (where $1 \le b < d$), noting that when $b = d$, then $\fM_4 = \{ 0\}$, in which case -- relative to the decomposition $\hilb = \fM_1 \oplus \fM_2 \oplus \fM_3$:
\[
J= \begin{bmatrix}
	Q_{1, 1} & Q_{1, 2} & 0 \\
	0 & 0 & I_\mu \\
	0 & 0 & 0 \end{bmatrix}, \]
while
\[
M =
	\begin{bmatrix}
	M_{1, 1} & M_{1, 2} & M_{1, 3} \\
	0 & M_{2, 2} & M_{2, 3} \\
	0 & 0 & M_{3, 3} \end{bmatrix},
	\ \ \ \ \ \ \ \ \ \
N =
	\begin{bmatrix}
	N_{1, 1} & N_{1, 2} & N_{1, 3} \\
	0 & N_{2, 2} & N_{2, 3} \\
	0 & 0 & N_{3, 3} \end{bmatrix}. \]
Arguing as before, this implies that $\begin{bmatrix} 0 & I_\mu \\ 0 & 0 \end{bmatrix}$ is a commutator of the two square-zero operators
\[
\begin{bmatrix} 		
	M_{2, 2} & M_{2, 3} \\
	0 & M_{3, 3} \end{bmatrix}, \ \ \ \
\text{ and }  \ \ \ \ \
\begin{bmatrix}
	N_{2, 2} & N_{2, 3} \\
	0 & N_{3, 3} \end{bmatrix}. \]
By Proposition~\ref{prop3.06}, this can happen only if $\mu$ is even.  In this case, however,  $\mu = \sum_{i=1}^d \mu_i$.

\bigskip

The result now follows by applying Theorem~\ref{thm3.08}, the reformulation of the Theorem of Borwein and Richmond.
\end{pf}

%%%%%%%%%%%%%%%%%%%%%%%%%%%%%%%%%%%%%%%%%%
%
%\begin{lem} \label{lem2.15}
%Suppose that $m \ge n$ and that $m-n \le 1$.   Then $J_m \oplus J_n \in \cniltwo$.
%\end{lem}
%
%\begin{pf}
%	Let $\{ e_k \}_{k=1}^{m+n}$ denote the standard orthonormal basis for $\bbC^{m+n}$.   Define $M \in \cB(\bbC^{m+n})$ as follows:
%	\[
%	M e_{j} = \begin{cases} 0 & \text { if } 1 \le j \le m+n \text{ is odd, and } \\ e_{j-1} & \text{ if } 1 \le j \le m+n  \text{ is even}, \end{cases}\]
%	while
%	\[
%	N e_{j} = \begin{cases} 0 & \text { if } 2 \le j \le m+n \text{ and } $j$ \text{ is even or } $j = 1$, \text{ and } \\ e_{j-1} & \text{ if } 3 \le j \le m+n  \text{ and } $j$ \text{ is odd}. \end{cases}\]
%	A routine computation shows that $[M, N] \sim J_m \oplus -J_n \sim J_m \oplus J_n$.    Thus $J_m \oplus J_n \in \cniltwo$.
%\end{pf}

%%%%%%%%%%%%%%%%%%%%%%%%%%%%%%%%%%%%%%%%%%

We next turn our attention to proving that the converse of Theorem~\ref{thm3.12} also holds.

\bigskip

\begin{prop} \label{prop3.13}
Let $\kappa \ge 1$.   Then $J_\kappa \oplus J_\kappa$ and $J_{\kappa+1} \oplus J_\kappa  \in \cniltwo$.
\end{prop}

\begin{pf}
Set $T_1 := J_\kappa \oplus J_\kappa$ and $T_2 := J_{\kappa+1} \oplus J_\kappa$.

In the case of $T_1$, simply note that $J_\kappa \sim -J_\kappa$, and thus $T_1 \sim J_\kappa \oplus -J_\kappa$.   Since $\cniltwo$ is invariant under similarity, we then have that
\[
T_1 \sim J_\kappa \oplus -J_\kappa  = \left[ \begin{bmatrix} 0 & I_\kappa \\ 0 & 0 \end{bmatrix}, \begin{bmatrix} 0 & 0 \\ J_\kappa & 0 \end{bmatrix} \right ] \in \cniltwo. \]

\bigskip

In the case of $T_2$, set $A = \begin{bmatrix} 0 \\ I_\kappa \end{bmatrix} \in \bbM_{\kappa+1, \kappa}(\bbC)$ and $Y = \begin{bmatrix} I_\kappa & 0 \end{bmatrix} \in \bbM_{\kappa, \kappa+1}(\bbC)$.  If we set
\[
M = \begin{bmatrix} 0 & A \\ 0 & 0 \end{bmatrix}, N = \begin{bmatrix} 0 & 0 \\ Y & 0 \end{bmatrix} \in \cB(\bbC^{\kappa+1} \oplus \bbC^\kappa) \simeq \bbM_{2\kappa+1}(\bbC), \]
then
\[
[M, N] = AY \oplus - YA = J_{\kappa+1}^* \oplus - J_\kappa^* \sim J_{\kappa+1} \oplus J_\kappa = T_2.\]
Since $M, N \in \niltwo$, $T_2 \in \cniltwo$.
\end{pf}

%%%%%%%%%%%%%%%%%%%%%%%%%%%%%%%%%%%%%%%%%%

\begin{cor} \label{cor3.14}
Let $J \in \bbM_n(\bbC)$ be a nilpotent operator.   Suppose that
\[
J \sim Q := J_{\kappa_1} \oplus J_{\kappa_2} \oplus \cdots \oplus J_{\kappa_d}, \]
where $\sum_{j=1}^d \kappa_j = n$ and $\kappa_1 \ge \kappa_2 \ge \cdots \ge \kappa_d$.

If, for $1 \le j \le \lfloor \frac{n}{2} \rfloor$, we have that $\kappa_{2j-1} - \kappa_{2j} \in \{ 0, 1\}$ \emph{(}and $\kappa_d = 1$ in the event that $d$ is odd\emph{)}, then $J \in \cniltwo$.
\end{cor}

\begin{pf}
Define
\[
Q_{0} := \oplus_{j=1}^{\lfloor \frac{n}{2} \rfloor} (J_{\kappa_{2j-1}} \oplus J_{\kappa_{2j}}). \]
By Proposition~\ref{prop3.13}, each $(J_{\kappa_{2j-1}} \oplus J_{\kappa_{2j}}) \in \cniltwo$, from which it easily follows that $Q_{0} \in \cniltwo$.   If $d$ is even, then $Q = Q_0$ and we are done.
If $d$ is odd, then $Q= Q_{0} \oplus J_1$, which is clearly in $\cniltwo$, as $Q_{0}$ and $J_1$ are in $\cniltwo$.

\smallskip

Finally, $J \in \cniltwo$ because $J \sim Q$, and $\cniltwo$ is invariant under similarity.
\end{pf}

%%%%%%%%%%%%%%%%%%%%%%%%%%%%%%%%%%%%%%%%%%

\begin{thm} \label{thm3.15}
Let $2 \le n \in \bbN$ and suppose that $J \in \cB(\bbC^n)$ is a nilpotent operator.   The following conditions are equivalent.
\begin{enumerate}
	\item[(a)]
	$J = [M, N]$ where $M^2 = 0 = N^2$; i.e. $J \in \cniltwo$.
	\item[(b)]
	$J$ admits a square root.
%	\item[(c)]
%	$J_+ := MN + NM$ admits a square root.
\end{enumerate}	
\end{thm}

\begin{pf}
\begin{enumerate}
	\item[(a)] implies (b).
	
	This is the content of Theorem~\ref{thm3.12}.
	\item[(b)] implies (a).
	
	Choose integers $1 \le d$ and $\kappa_1 \ge \kappa_2 \ge \cdots \ge \kappa_d$ such that $J = \oplus_{i=1}^d J_{\kappa_i}$.
	Suppose that $J$ admits a square root.  By Theorem~\ref{thm3.07}, $\kappa_{2 j -1} - \kappa_{2 j} \le 1$ for all $1 \le j \le \lfloor \frac{d}{2} \rfloor$, and if $d$ is odd, then $\kappa_d = 1$, implying that $J_{\kappa_d} = 0$.
	
	\smallskip
	
	By Corollary~\ref{cor3.14}, $J \in \cniltwo$.
\end{enumerate}	
\end{pf}

\smallskip

As an immediate consequence of the above Theorem, we retrieve the result of Corollary~\ref{cor3.04}  that if an element of $\cniltwo \cap \cB(\bbC^n)$ is nilpotent, then it is nilpotent of order at most $\lfloor \frac{n+1}{2} \rfloor$.   The question of which nilpotent (or general) operators admit $m^{th}$ roots is an interesting one.   We direct the reader to~\cite{Otero1990, Ozturk2021, Psarrakos2002, Yood2002} for more information on this topic.

%%%%%%%%%%%%%%%%%%%%%%%%%%%%%%%%%%%%%%%%%%

\begin{cor} \label{cor3.16}
Let $2 \le n \in \bbN$ and let $T \in \cB(\bbC^n)$.   If $T \in \cniltwo$, then  $T$ admits a square root.
\end{cor}

\begin{pf}

	Since both $\cniltwo$ and the set $\{ Z^2 : Z \in \cB(\bbC^n \}$ are invariant under similarity, we see from Remark~\ref{rem1.02}  that without loss of generality, we may assume that $T = B \oplus Q$, where $B$ is invertible and $Q$ is nilpotent.    By Proposition~\ref{prop2.02}, the hypothesis that $T \in \cniltwo$ implies that $B \in \cniltwo$ and $Q \in \cniltwo$.
	
	Any invertible operator $B \in \cB(\hilb)$ (where $\dim\, \hilb < \infty$) admits a square root $X$;  one can either appeal to Theorem~1 of~\cite{BorweinRichmond1984}, or one may simply apply the functional calculus to $B$ using a branch of the square-root function $f(z) = z^{\frac{1}{2}}$ that avoids the finitely many non-zero eigenvalues of $B$.   Theorem~\ref{thm3.15} above shows that if $Q \in \cniltwo$ is nilpotent and acts on a finite-dimensional space, then $Q = Y^2$ for some nilpotent operator $Y$.
	
	Thus $T = B \oplus Q = X^2 \oplus Y^2 = (X \oplus Y)^2$.
\end{pf}

Note that the converse of this result is false.   The identity operator $I \in \cB(\bbC^n)$ obviously admits a square root, but it is not a commutator (of any two operators) as its trace is non-zero.

%%%%%%%%%%%%%%%%%%%%%%%%%%%%%%%%%%%%%%%%%%

\begin{thm} \label{thm3.17}
Let $2 \le n \in \bbN$ and $T \in \cB(\bbC^n)$.   The following are equivalent.
\begin{enumerate}
	\item[(a)]
	$T \in \cniltwo$;
	\item[(b)]
	There exists $0 \le 2 k \le n$, $A \in \cB(\bbC^k)$ invertible and $L \in \cB(\bbC^{n - 2k})$ nilpotent such that
	\[
	T \sim (A \oplus - A) \oplus L^2. \]
\end{enumerate}
\end{thm}	

\begin{pf}
\begin{enumerate}
	\item[(a)] implies (b). \ \ \
	Suppose that $T \in \cniltwo$.  By Remark~\ref{rem1.02}, we may assume without loss of generality that $T = B \oplus Q$, where $B$ is invertible and $Q$ is nilpotent.    By Proposition~\ref{prop2.02}, $B \in \cniltwo$ and $Q \in \cniltwo$.  (One of $B$ and $Q$ may be absent.  We deal with the case where both are present.)   If $B$ is present, then by Proposition~\ref{prop2.03}, $B$ acts on a space of even dimension, say $2k$, and there exists $A \in \cB(\bbC^k)$ such that $B \sim A \oplus -A$.
	By Theorem~\ref{thm3.12}, $Q = L^2$ for some nilpotent operator $L \in \cB(\bbC^{n-2k})$.
	From this the result immediately follows.	
	\item[(b)] implies (a). \ \ \
	If $T \sim (A \oplus -A) \oplus L^2$ where $A$ is invertible and $L$ is nilpotent, then setting $B := (A \oplus -A)$ and $Q = L^2$, we see from Proposition~\ref{prop2.03} and Theorem~\ref{thm3.12} that $B$ and $Q$ lie in $\cniltwo$, whence $T \in \cniltwo$.
\end{enumerate}
\end{pf}

%%%%%%%%%%%%%%%%%%%%%%%%%%%%%%%%%%%%%%%%%%

It is interesting to compare the above theorem with the following result of Smith~\cite[Theorem~1]{Smith1976}.

\begin{thm}  \emph{(\textbf{Smith}) } \label{thm3.18}
Let $n\geq 3$, $T\in \bbM_n(\bbC)$ with $\textsc{Tr}(T)=0$. Then there are nilpotent
matrices $A, B\in \bbM_n(\bbC)$ such that $C=AB-BA$.
\end{thm}

%%%%%%%%%%%%%%%%%%%%%%%%%%%%%%%%%%%%%%%%%%

\begin{rem} \label{rem2.24}
It is interesting to note that if $M, N \in \cB(\bbC^n)$ are square-zero operators, and if $J := [M, N]$ is nilpotent, then by Theorem~\ref{thmRR2002}, $M$ and $N$ are simultaneously triangularisable.
But then the Jordan product $J_+ := MN + N M = (M+N)^2$ is nilpotent, and so by Theorem~\ref{thm3.17}, $J_+ \in \cniltwo$ as well.
\end{rem}

%%%%%%%%%%%%%%%%%%%%%%%%%%%%%%%%%%%%%%%%%%

%%%%%%%%%%%%%%%%%%%%%%%%%%%%%%%%%%%%%%%%%%
%%%%%%%%%%%%%%%%%%%%%%%%%%%%%%%%%%%%%%%%%%
%%%%%%%%%%%%%%%%%%%%%%%%%%%%%%%%%%%%%%%%%%

%%%%%%%%%%%

%%%%%%%%%%%%%%%%%%%%%%%%%%%%%%%%%%%%%%%%%%
%%%%%%%%%%%%%%%%%%%%%%%%%%%%%%%%%%%%%%%%%%
% SECTION FOUR
%%%%%%%%%%%%%%%%%%%%%%%%%%%%%%%%%%%%%%%%%%
%%%%%%%%%%%%%%%%%%%%%%%%%%%%%%%%%%%%%%%%%%

\vskip 1.5 cm

\section{The closure of $\cniltwo$} \label{section4}

%%%%%%%%

%%%%%%%%%%%%%%%%%%%%%%%%%%%%%%%%%%%%%%%%%%
%%%%%%%%%%%%%%%%%%%%%%%%%%%%%%%%%%%%%%%%%%
%%%%%%%%%%%%%%%%%%%%%%%%%%%%%%%%%%%%%%%%%%

%%%%%%%%%%%%%%%%%%%%%%%%%%%%%%%%%%%%%%%%%%

\subsection{} \label{sec4.01}
Determining a complete list of invariants for a similarity-invariant subset $\Gamma$ of $\bofh$ when $\hilb$ is infinite-dimensional is, in general, an extremely difficult problem.   The problem is much more tractable if the similarity-invariant set is norm-closed.  Indeed, in the 1980s, Herrero~\cite{Herrero1989} and  Apostol, Fialkow, Herrero and Voiculescu~\cite{ApostolFialkowHerreroVoiculescu1984} developed an extensive machinery to study the \emph{norm-closures} of similarity-invariant subsets of $\bofh$.    Herrero~\cite{Herrero1990} established a meta-theorem which (\emph{very} roughly) states that to calculate the norm-closure of a similarity-invariant subset $\Gamma$ of $\bofh$, one enumerates the  spectral, semi-Fredholm index and algebraic conditions necessary for an operator to belong to $\ttt{clos}(\Gamma)$, and in most cases these conditions should also be sufficient.

Following in the footsteps of Apostol et al, we turn our attention to characterising the norm-closure of $\cniltwo$.   In the finite-dimensional setting, we obtain a complete characterisation of this set, while in the infinite-dimensional setting, we determine the intersection of $\ttt{clos}(\cniltwo)$ with the set $\ttt{(bqt)}$ of biquasitriangular operators.  %We shall also exhibit an index obstruction to belonging to $\ttt{clos}(\cniltwo)$, and we use this to demonstrate that very few weighted shift operators belong to that set.

\bigskip
%
%We begin with the finite-dimensional case, where we shall require a result of Horn and Johnson~\cite[Theorem~1.3.20]{HornJohnson1990}.

%%%%%%%%%%%%%%%%%%%%%%%%%%%%%%%%%%%%%%%%%%

%\begin{thm} \emph{\textbf{(Horn and Johnson)}} \ \ \ \label{thm4.02}
%
%Let $\hilb, \cK$ be Hilbert spaces, $A \in \bofh$ and $B \in \cB(\cK)$.   Then $AB$ and $BA$ have the same non-zero eigenvalues, counted according to algebraic multiplicity.   If $\hilb = \cK$ and either $A$ or $B$ is invertible, then $AB$ and $BA$ are similar.
%\end{thm}

%
%\begin{pf}
%Observe that
%\[
%	\begin{bmatrix} AB & 0 \\ B & 0 \end{bmatrix} \begin{bmatrix} I & A \\ 0 & I \end{bmatrix} = \begin{bmatrix} A B & ABA \\ B & BA \end{bmatrix} = \begin{bmatrix} I & A \\ 0 & I \end{bmatrix} \begin{bmatrix} 0 & 0 \\ B & BA \end{bmatrix}. \]
%
%Since $ \begin{bmatrix} I & A \\ 0 & I \end{bmatrix}$ is invertible, we see that $C := \begin{bmatrix} AB & 0 \\ B & 0 \end{bmatrix}$ and $D := \begin{bmatrix} 0 & 0 \\ B & BA \end{bmatrix}$ are similar.  In particular, they have the same eigenvalues, repeated with the same algebraic (and geometric) multiplicities.
%
%But the non-zero eigenvalues of $C$ are precisely the non-zero eigenvalues of $AB$  (counted according to algebraic  multiplicity), while the non-zero eigenvalues of $D$ are precisely the non-zero eigenvalues of $BA$ (counted according to algebraic multiplicity), from which the result follows.
%
%\smallskip
%
%The last assertion follows from the fact that (if we call the invertible operator $A$, then) $A B = A (BA) A^{-1}$.
%\end{pf}

%%%%%%%%%%%%%%%%%%%%%%%%%%%%%%%%%%%%%%%%%%

\begin{thm} \label{thm4.03}
Suppose that $\dim\, \hilb < \infty$.   Then
\[
\ttt{clos} (\mathfrak{c} (\ttt{nil}_2)) = \ttt{Bal}(\hilb). \]
\end{thm}

\begin{pf}
Suppose that $T = [M, N]$, where $M^2 = 0 = N^2$.   By Theorem~\ref{thm3.17}, there exist an invertible operator $A$ and a nilpotent operator $L$ such that
\[
T \sim (A \oplus -A) \oplus L^2. \]
It follows directly from this that $T \in \ttt{Bal}(\hilb)$.   That is, $\cniltwo \subseteq \ttt{Bal}(\hilb)$.

Since $\ttt{Bal}(\hilb)$ is closed \cite[Theorem~3.5]{MarcouxRadjaviZhang2024},
\[
\ttt{clos}(\cniltwo) \subseteq \ttt{Bal}(\hilb). \]

\smallskip

Conversely, suppose that $T \in \ttt{Bal}(\hilb)$, and let $\eps > 0$.   The proof of Theorem~3.5 of~\cite{MarcouxRadjaviZhang2024} shows we can find an operator $X \in \bofh$ such that $\norm X - T \norm < \eps$, and $X$ is similar to a diagonal operator $D = \ttt{Diag}(\beta_1, \beta_2, \ldots, \beta_n)$ such that for all $1 \le k \le \lfloor\frac{n}{2}\rfloor$, $\beta_{2k} = -\beta_{2k-1}$, $\beta_n = 0$ if $n$ is odd, and for $1 \le i, j \le n$, $\beta_i = \beta_j$ if and only if $i = j$.
In particular,
\[
D \simeq \oplus_{k=1}^{\lfloor n/2 \rfloor} \ttt{Diag}(\beta_{2k-1}, -\beta_{2k-1}) \oplus 0^{(m)}, \]
where $m = \begin{cases} 0 & \text{ if } n \text{ is even } \\ 1 & \text{ if } n \text{ is odd.}\end{cases}$

Note that $\begin{bmatrix} \alpha & 0 \\ 0 & - \alpha \end{bmatrix} \in \cniltwo$ for  each $\alpha \in \bbC$, from which we  conclude that $D \in \cniltwo$.  Since $\cniltwo$ is closed under similarity, $X \in \cniltwo$, whence $T \in \ttt{clos}(\cniltwo)$.   That is,
\[
\ttt{clos}(\cniltwo) \supseteq \ttt{Bal}(\hilb). \]

This completes the proof of the Proposition.
\end{pf}

%%%%%%%%%%%%%%%%%%%%%%%%%%%%%%%%%%%%%%%%%%

\smallskip

The set $\ttt{Neg}_S(\hilb) := \{ T \in \bofh: T \sim -T \}$  was originally defined in~\cite{MarcouxRadjaviZhang2024}, where it was shown to contain the set $\fC_\fE := \{ [E, F] : E, F \in \bofh, E^2 = E, F^2 = F \}$.

\begin{cor}  \label{cor4.04}
If $\dim\, \hilb < \infty$, then
\[
\ttt{clos}(\cniltwo) = \ttt{clos}(\fC_\fE) = \ttt{clos}(\ttt{Neg}_{S}(\hilb)) = \ttt{Bal}(\hilb). \]
\end{cor}

\begin{pf}
This follows immediately from Corollary~3.6 of~\cite{MarcouxRadjaviZhang2024}.
\end{pf}

We remark in passing that every nilpotent operator in $\cB(\bbC^n)$ is balanced, and thus lies in $\ttt{clos}(\cniltwo)$.

%%%%%%%%%%%%%%%%%%%%%%%%%%%%%%%%%%%%%%%%%%

\subsection{} $\ttt{clos}\, (\cniltwo)$ \textbf{in the infinite-dimensional setting.} \label{sec4.05}
As indicated above, a complete characterisation of which Hilbert space operators lie in $\cniltwo$ when $\hilb$ is infinite-dimensional seems out of reach at the moment.    We begin with an example that shows that  $\cniltwo$  fails to be invariant under approximate unitary equivalence.

%%%%%%%%%%%%%%%%%%%%%%%%%%%%%%%%%%%%%%%%%%

\begin{eg} \label{eg4.06}
There exist normal operators $L_1$ and $L_2 \in \bofh$ such that $L_1 \simeq_a L_2$, $L_1 \in \cniltwo$, but $L_2 \not \in \cniltwo$.

Let $D = \ttt{diag} (d_n)_n \in \bofh$ be a diagonal operator with $d_1 = 1$ and $\{ d_n\}_{n=2}^\infty$ dense in $[\frac{1}{2}, 1)$. Then $\sigma(D) = [\frac{1}{2}, 1]$ and $1$ is an eigenvalue of $D$ of (geometric and algebraic) multiplicity $1$.    Let $M_x$ denote the usual multiplication operator on $L^2([\frac{1}{2}, 1], dm)$, where $m$ denotes Lebesgue measure.

Let $L_1 := (-D) \oplus D$ and $L_2 = (-D) \oplus M_x$.    Since $\sigma(L_1) = [-1, -\frac{1}{2}] \cup [\frac{1}{2}, 1] = \sigma(L_2)$, it follows from the Weyl-von Neumann-Berg Theorem~\cite[Theorem~II.4.4]{Davidson1996} that $L_1 \simeq_a L_2$.  Note that $\sigma_p(L_2) = \sigma_p(-D) \subseteq  [-1, -\frac{1}{2}]$, and so clearly $L_2$ is not similar to $-L_2$.

By Theorem~\ref{thm2.06}, $L_1 \in \cniltwo$, while  $L_2 \not \in \cniltwo$.
	
\end{eg}

\subsection{} \textbf{Biquastriangular operators in}  $\ttt{clos}({\mathfrak{c}(\ttt{nil}_2)})$\textbf{.} \label{sec4.07}

Let us next consider which biquastriangular operators lie in $\ttt{clos}({\mathfrak{c}(\ttt{nil}_2)})$.

%%%%%%%%%%%%%%%%%%%%%%%%%%%%%%%%%%%%%%%%%%
%%%%%%%%%%%%%%%%%%%%%%%%%%%%%%%%%%%%%%%%%%

\smallskip

The following is a generalisation of Lemma 3.2 of the paper~\cite{MarcouxRadjaviZhang2023}.

\begin{prop}\label{prop4.08}
Let $\cA$ be a unital Banach algebra.   Let $a, b_n \in \cA$, $n \ge 1$, and suppose that $a = \lim_n b_n$.   Suppose that for each $n \ge 1$,  every connected component of $\sigma(b_n)$ intersects $\sigma(-b_n)$.
Then every connected component of $\sigma(a)$ intersects $\sigma(-a)$.
\end{prop}

\begin{pf}
Suppose to the contrary that there exists a connected component $\Omega \subseteq \sigma(a)$ satisfying
	\[
	\Omega \cap -\sigma(a) = \Omega \cap \sigma(-a) = \varnothing.   \]
Note that $\Omega$ and $\sigma(-a)\cup(\sigma(a)\setminus \Omega)$ are  disjoint, compact sets of $\bbC$.   Thus, we can find  disjoint open sets $G_1, G_2$ of $\bbC$ such that $\Omega \subseteq G_1$, $\sigma(-a)\cup(\sigma(a)\setminus \Omega) \subseteq G_2$ and $\textup{d}(G_1,G_2)=\delta>0$.
		
	By applying the upper semi-continuity of the spectrum and Newburgh's Theorem (c.f.~\cite[Theorem~1.1]{Herrero1989}), we find that there exists $\delta > 0$ such that $\norm y  - a \norm < \delta$ implies that
\[\sigma(y)\subset G_1\cup G_2, \sigma(-y) \subseteq G_2, \textup{and}~\sigma(y) \cap G_1 \ne \varnothing.\]  In particular,  for $N \in \bbN$ sufficiently large we have that
\[\sigma(b_N)\subset G_1\cup G_2, \sigma(-b_N) \subseteq G_2, \textup{and}~\sigma(b_N) \cap G_1 \ne \varnothing.\]
%Of course, a connected component of $\sigma(b_N)$ is either contained in $G_1$ or it is contained in $G_2$.
Let  $\beta \in \sigma(b_N) \cap G_1 \ne \varnothing$.  Clearly the connected component $\Gamma_\beta$ of $\beta$ in $\sigma(b_N)$ also lies in $G_1$.  But then
\[
\Gamma_\beta \cap \sigma(-b_N) \subseteq G_1 \cap G_2 = \varnothing, \]
contradicting our hypothesis.  This completes the proof.
%But then any component of $\sigma(b_N)$ which lies in $G_1$ has no intersection with $G_2$, and hence with $\sigma(-b_N)$. However, $\sigma(b_N) \cap G_1 \ne \varnothing$ implies that there is indeed a component contained in $G_1$. This contradicts with the assumption.
\end{pf}

%%%%%%%%%%%%%%%%%%%%%%%%%%%%%%%%%%%%%%%%%%

\bigskip

By combining Proposition~\ref{prop2.12} and Proposition \ref{prop4.08}, one obtains

\begin{prop}\label{prop4.09}
Let $T \in \bofh$. Suppose that $T\in \ttt{clos}({\mathfrak{c}(\ttt{nil}_2)})$. Then
\begin{enumerate}
	\item[(a)]
	every connected component of $\sigma(T)$ intersects $\sigma(-T)$;
	\item[(b)] every connected component of $\sigma_e(T)$ intersects $\sigma_e(-T)$.
\end{enumerate}
\end{prop}

%%%%%%%%%%%%%%%%%%%%%%%%%%%%%%%%%%%%%%%%%%

\begin{prop}\label{prop4.10}
Let $T \in \bofh$. Suppose that $T\in \ttt{clos}({\mathfrak{c}(\ttt{nil}_2)})$. If $\alpha, -\alpha$ are isolated points,
then
\[
		\dim\, \hilb(\{ \alpha\}; T) = \dim\, \hilb(\{ -\alpha\}; T). \]
\end{prop}

\begin{pf}
Without loss of generality, we may assume that $\alpha\neq 0$ and that
\[\min\{ \dim\, \hilb(\{ \alpha\}; T), \dim\, \hilb(\{ -\alpha\}; T)\}<+\infty.\]
 We will consider the following two subcases.
\begin{enumerate}
	\item[(a)]
	If $\dim\, \hilb(\{ \alpha\}; T)=+\infty$, then $\alpha\in \sigma_e(T)$, but $-\alpha\notin \sigma_e(T)$ and hence $\alpha\notin \sigma_e(-T)$. Then $\{\alpha\}$, which is a component of $\sigma_e(T)$, has no intersection with $\sigma_e(-T)$. This, however, contradicts Proposition \ref{prop4.09}.  A similar argument may be applied  if\linebreak $\dim\, \hilb(\{-\alpha\}; T) = +\infty$.
	\item[(b)]
	Next, suppose that
\[\max\{ \dim\, \hilb(\{ \alpha\}; T), \dim\, \hilb(\{ -\alpha\}; T)\}<+\infty.\]
Choose $r\in (0, \frac{|\alpha|}{3})$, set $G_1=B(\alpha,r) := \{ \lambda \in \bbC: |\lambda - \alpha| < r\}$, and $G_2=B(-\alpha,r)$.
Then $G_1,G_2$ are disjoint open sets.  By Corollary 1.6 of \cite{Herrero1989},
there exists $\delta>0$ such that $\|Y-T\|<\delta$ implies that
\begin{itemize}
	\item{}
	$\sigma(Y) \cap G_k \ne \varnothing$, $k = 1, 2$,
	\item{}
	$\dim\, \hilb(\sigma(Y) \cap G_1; Y)= \dim\, \hilb(\{ \alpha\}; T)<\infty$, and
    \item{}
 $\dim\, \hilb(\sigma(Y) \cap G_2; Y) =\dim\,\hilb(\{ - \alpha\}; T)<\infty$.
\end{itemize}
Since $T\in \ttt{clos}({\mathfrak{c}(\ttt{nil}_2)})$, we can choose $X\in {\mathfrak{c}(\ttt{nil}_2)}$, such that $\|X-T\|<\delta$. By Theorem~\ref{thm2.08},
\[\dim\, \hilb(\sigma(X) \cap G_1; X)=\dim\, \hilb(\sigma(X) \cap G_2; X).\]
Therefore, \[
		\dim\, \hilb(\{ \alpha\}; T) = \dim\, \hilb(\{ -\alpha\}; T). \]
\end{enumerate}
\end{pf}

%%%%%%%%%%%%%%%%%%%%%%%%%%%%%%%%%%%%%%%%%%

\begin{thm}\label{thm4.11}
Let $T \in \bofh$ be a biquasitriangular operator.   The following conditions are equivalent:
\begin{enumerate}
	\item[(a)]
	$T \in \ttt{clos}\, (\fC_\fE)$;
	\item[(b)]
	$T \in \ttt{clos}\, (\ttt{Bal}(\hilb))$;
	\item[(c)]
	\begin{enumerate}
		\item[(i)]
		every connected component of $\sigma(T)$ intersects $\sigma(-T)$;
		\item[(ii)]
		every connected component of $\sigma_e(T)$ intersects $\sigma_e(-T)$; and
		\item[(iii)]
		If $\alpha, -\alpha \in \sigma(T)$ are isolated points, then
		\[
		\dim\, \hilb(\{ \alpha\}; T) = \dim\, \hilb(\{ -\alpha\}; T); \]
	\end{enumerate}
\item [(d)] $T \in \ttt{clos}({\mathfrak{c}(\ttt{nil}_2)})$.
\end{enumerate}	

\end{thm}

\begin{pf}
By Theorem 4.11 of~\cite{MarcouxRadjaviZhang2023} (a), (b), and (c) are equivalent. By Proposition \ref{prop4.09} and Proposition \ref{prop4.10}, (d) implies (c).

 By examining  the proof of (c) implying (a) of Theorem~3.5 and Theorem~4.11 of ~\cite{MarcouxRadjaviZhang2023}, we find that to show that (c) implies (d), it suffices to note that balanced operators acting on a finite dimensional Hilbert space are in $\ttt{clos}({\mathfrak{c}(\ttt{nil}_2)})$(Corollary \ref{cor4.04})) and to show that $R := X\oplus -X\in {\mathfrak{c}(\ttt{nil}_2)}$.
 The last fact is easy to check by factoring $X\oplus -X$
 as
 \[\begin{bmatrix}X&\\& -X \end{bmatrix}=\begin{bmatrix} 0& X\\ & 0 \end{bmatrix}\begin{bmatrix} 0& \\I & 0 \end{bmatrix}-\begin{bmatrix} 0& \\I & 0 \end{bmatrix}\begin{bmatrix} 0& X\\ & 0 \end{bmatrix}.\]
\end{pf}

%%%%%%%%%%%%%%%%%%%%%%%%%%%%%%%%%%%%%%%%%%

\begin{cor} \label{cor4.12}
Let $\hilb$ be an infinite-dimensional, separable Hilbert space.  Then\linebreak $\ttt{clos}({\mathfrak{c}(\ttt{nil}_2)})$ contains the closure of the set of all nilpotent operators in $\bofh$.
\end{cor}

%%%%%%%%%%%%%%%%%%%%%%%%%%%%%%%%%%%%%%%%%%
%%%%%%%%%%%%%%%%%%%%%%%%%%%%%%%%%%%%%%%%%%

%\subsection{} \label{sec3.08}
%Our goal in this section is to characterise the following sets.
%\begin{enumerate}
%	\item[(a)]
%	$\kofh \cap \ttt{clos}(\cniltwo)$;
%	\item[(b)]
%	$\ttt{(bqt)} \cap \ttt{clos}(\cniltwo)$; and
%	\item[(c)]
%	$\ttt{(Nor)} \cap \ttt{clos}(\cniltwo)$, where $\ttt{(Nor)} := \{ L \in \bofh:  L \text{ is normal}\}$.
%\end{enumerate}
%Intriguingly, we shall demonstrate that
%\[
%\ttt{(bqt)} \cap \ttt{clos}(\cniltwo) = \ttt{(bqt)} \cap \ttt{clos}(\fE) = \ttt{(bqt)} \cap \ttt{clos}(\ttt{Bal}(\hilb)). \]
%
%Furthermore, we exhibit an index obstruction to belonging to $\ttt{clos}(\cniltwo)$, and we use this to show that very few weighted shift operators belong to this set.
%
%%%%%%%%%%%%%%%%%%%%%%%%%%%%%%%%%%%%%%%%%%%
%%%%%%%%%%%%%%%%%%%%%%%%%%%%%%%%%%%%%%%%%%

\subsection{}\ \ \  \label{sec4.13}
We next turn our attention to determining which compact operators lie in the set $\ttt{clos}(\cniltwo)$.   Interestingly, we shall show that a compact operator is a limit of commutators of square-zero operators if and only if it is a limit of \emph{compact} square-zero operators.   The corresponding statement where the adjective ``compact" is replaced by ``normal" turns out to be false (see Example~\ref{eg4.16} below).

%%%%%%%%%%%%%%%%%%%%%%%%%%%%%%%%%%%%%%%%%%

\begin{lem} \label{lem4.14}
Let $K \in \kofh$.   If $K \in  \ttt{clos}(\cniltwo)$, then $K$ is balanced.
\end{lem}

\begin{pf}
This is an immediate consequence of the equivalence of (c) and (d) in Theorem~\ref{thm4.11}.
\end{pf}

%%%%%%%%%%%%%%%%%%%%%%%%%%%%%%%%%%%%%%%%%%

\begin{thm} \label{thm4.15}
Let $K \in \kofh$.   The following are equivalent.
\begin{enumerate}
	\item[(a)]
	$K \in \ttt{clos}(\cniltwo)$.
	\item[(b)]
	$K \in \ttt{clos}(\ttt{Bal}(\hilb))$.
	\item[(c)]
	$K \in \ttt{Bal}(\hilb)$.
	\item[(d)]
	$K \in \ttt{WBal}(\hilb)$.
	\item[(e)]
	$K \in \ttt{clos}(\kofh \cap \cniltwo)$.
\end{enumerate}	
\end{thm}

\begin{pf}
\begin{enumerate}
	\item[(a)] if and only if  (b). \ \ \ This follows immediately from Theorem~\ref{thm4.11} above.
	\item[(b)] implies (c). \ \ \ This is Lemma~\ref{lem4.14}.
	\item[(c)] implies (d). \ \ \ This follows from the definitions of both sets.
	\item[(d)] implies (e). \ \ \ Let $K \in \ttt{WBal}(\hilb)$ be a compact operator.   We shall show that given $\eps> 0$, there exists a finite-dimensional subspace $\cM \subseteq \hilb$ and $F \in \ttt{Bal}(\cM)$ such that $\norm K - (F \oplus 0) \norm < \frac{\eps}{2}$.   If we can do this, then by applying Theorem~\ref{thm4.03}, we can find $B \in \cniltwo \cap \cB(\cM)$ such that $\norm B - F \norm< \frac{\eps}{2}$, whence $B \oplus 0 \in \cniltwo \cap \kofh$, and $\norm K - (B \oplus 0) \norm < \eps$.   This implies that $\ttt{WBal}(\hilb) \cap \kofh \subseteq \ttt{clos}(\kofh \cap \cniltwo)$, as required.
	
	\bigskip

	The spectrum of any compact operator in $\bofh$ is either finite, or a sequence converging to $0$.   Suppose that $K \in \ttt{WBal}(\hilb)$.   It easily follows that we must have $\sigma(K) = \{ \beta_i, -\beta_i\}_{i=1}^\eta \cup \{ 0\}$, where $\eta \in \bbN \cup \{ 0, \infty\}$, $\beta_i \ne 0$ if $1 \le i < \eta+1$,  and $\{ \beta_i, -\beta_i\} \cap \{ \beta_j, -\beta_j\} = \varnothing$ if $1 \le i \ne j < \eta+1$.   Moreover, if $\eta = \infty$, then $\lim_i \beta_i = 0$, and in all cases, $\mu(\beta_i) = \mu(-\beta_i)$ for all $1 \le i < \eta + 1$.
	
	Let $\eps > 0$.   Fix $0 \le p < \eta + 1$ such that $| \beta_i| \ge \frac{\eps}{4}$ if $1 \le i \le p$ and $| \beta_i | < \frac{\eps}{4}$ if $p < i < \eta + 1$.
	
	By considering the Riesz subspaces $\hilb_i$ corresponding to the clopen sets $\{ \beta_i, -\beta_i\}$, $1 \le i \le p$, we may decompose $\hilb = \hilb_1 \oplus \hilb_2 \oplus \cdots \oplus \hilb_p \oplus \hilb_0$ such that with respect to this decomposition,
	\[
	K = \begin{bmatrix} K_{11} & K_{12} & K_{13} & \cdots & K_{1p} & K_{1 0} \\
		& K_{22} & K_{23} & \cdots & K_{2p} & K_{2 0} \\
		& & \ddots &  & & \vdots \\
		& & & \ddots & \vdots & \vdots \\
		& & & & K_{p p} & K_{p 0} \\
		& & & & & K_{0 0}
		\end{bmatrix}, \]
	where $\sigma(K_{ii}) = \{ \beta_i, -\beta_i\}$, $1 \le i \le p$ and $\sigma(K_{00}) = \sigma(K) \setminus \{ \beta_j, -\beta_j\}_{j=1}^p \subseteq \bbD_{\eps/4} := \{ z \in \bbC: |z| < \frac{\eps}{4} \}$.

	Let $P$ denote the orthogonal projection of $\hilb$ onto $\oplus_{j=1}^p \hilb_j$, and observe that $\rank\, P < \infty$.   Since $K$ and $K_{00}$ are compact, we can find a finite-rank orthogonal projection $Q = (I-P)Q(I-P)$ such that
	\begin{itemize}
		\item{}
		$\norm (P+Q) K (P+Q) - K \norm < \frac{\eps}{4}$, and
		\item{}
		by furthermore appealing to the upper-semicontinuity of the spectrum,
		\[
		\sigma(Q K_{00} Q) \subseteq \bbD_{\eps/4}. \]
	\end{itemize}	
	(To achieve the second step, choose $\delta > 0$ such that $\norm X - K_{00} \norm < \delta$ implies that $\sigma(X) \subseteq \bbD_{\eps/4}$, and then choose $Q$ sufficiently ``large" such that $\norm Q K_{00} Q - K_{00}  \norm < \delta$.)
	
	Write $\cJ_0 = \ran\, Q$ and let $\cJ_1 := \hilb_0 \ominus \cJ_0$.  Relative to the decomposition $\hilb = \oplus_{j=1}^p \hilb_j \oplus \cJ_0 \oplus \cJ_1$, we have
	\[
	(P+Q) K (P + Q) =
	\begin{bmatrix} K_{11} & K_{12} & K_{13} & \cdots & K_{1p} & L_{1 0} & 0 \\
		& K_{22} & K_{23} & \cdots & K_{2p} & L_{2 0} & 0  \\
		& & \ddots &  & & \vdots & \vdots \\
		& & & \ddots & \vdots & \vdots & \vdots \\
		& & & & K_{p p} & L_{p 0} & 0  \\
		& & & & & L_{0 0} & 0 \\
		& & & & & 0 & 0
	\end{bmatrix}, \]
	where $\sigma(L_{00} \oplus 0) = \sigma(Q K_{00} Q)  \subseteq \bbD_{\eps/4}$.  Since $\dim\, \cJ_0 = \rank \, Q < \infty$, we can upper-triangularise $L_{00}$ with respect to some orthonormal basis for $\cJ_0$.   Let $D$ denote the diagonal of $L_{00}$, and observe that $\norm D \norm < \frac{\eps}{4}$, since the norm of a diagonal matrix is equal to its spectral radius.  Set $L := L_{00} - D$, so that $\sigma(L) = \{ 0 \}$, and set
	\[
	F :=    	
	\begin{bmatrix} K_{11} & K_{12} & K_{13} & \cdots & K_{1p} & L_{1 0} \\
		& K_{22} & K_{23} & \cdots & K_{2p} & L_{2 0}  \\
		& & \ddots &  & & \vdots  \\
		& & & \ddots & \vdots & \vdots \\
		& & & & K_{p p} & L_{p 0}  \\
		& & & & & L
	\end{bmatrix} \in \cB(\cM),  \]
	where $\cM = \oplus_{j=1}^p \hilb_j \oplus \cJ_0$.
	Then $\sigma(F) = \cup_{j=1}^p \sigma(K_{jj}) \cup \sigma(L)$, noting that the eigenvalues of $K_{jj}$ appear with the same algebraic multiplicity in $\sigma(F)$.  In other words, $F$ is balanced.  Finally,
	\begin{align*}
	\ \ \ \ \ \norm K - (F\oplus 0) \norm
		&\le \norm K - (P+Q) K (P+Q) \norm + \norm (P+Q) K (P+Q) - (F \oplus 0) \norm  \\
		&< \frac{\eps}{4} + \norm L - L_{00} \norm \\
		&< \frac{\eps}{4} + \frac{\eps}{4} \\
		&= \frac{\eps}{2}.
	\end{align*}
	As seen above, this completes the proof.	
	\item[(e)] implies (a). \ \ \ This is trivial.
	
\end{enumerate}
\end{pf}

%%%%%%%%%%%%%%%%%%%%%%%%%%%%%%%%%%%%%%%%%%

\begin{eg} \label{eg4.16}\ \ \ Let us denote the set of normal operators in $\bofh$ by $\ttt{(nor)}$.   After comparing conditions (a) and (e) of Theorem~\ref{thm4.15}, one might be tempted to ask whether or not $\ttt{clos} (\ttt{(nor)} \cap \cniltwo) = (\ttt{nor}) \cap \ttt{clos} (\cniltwo)$.   This is not the case, as we shall now see.

\bigskip

Suppose that $T=\begin{bmatrix} -I& \\ & D \\ \end{bmatrix}$, where $I\in \cB(\cH)$ is the identity operator and $D\in \cB(\cH)$ is
a normal operator with $\sigma_e(D)=\sigma(D)=[1,2]$. Then $T$ is normal and $T\in \clos({\mathfrak{c}(\ttt{nil}_2)})$ by Theorem~\ref{thm4.11}.
We claim that $T\notin \clos((\textsc{nor})\cap{\mathfrak{c}(\ttt{nil}_2)})$.

Otherwise, there exists a sequence $(T_n)_n$ of normal operators in ${\mathfrak{c}(\ttt{nil}_2)}$, such that $\underset{n\to \infty}{\lim}T_n=T$.
By Theorem~\ref{thm2.06}, each $T_n$ is unitarily equivalent to $-T_n$, hence in particular, $\sigma(T_n)=\sigma(-T_n)$.
By a result of Newburgh, the map taking an operator to its spectrum is continuous on the set of normal operators with the operator norm, to the set of compact subsets of $\bbC$ with the Hausdorff metric (see also Problem 105 of~\cite{Halmos1982}).
Therefore,
\[\sigma(T)=\underset{n\to \infty}{\lim}\sigma(T_n)=\underset{n\to \infty}{\lim}\sigma(-T_n)=\sigma(-T)=-\sigma(T),\]
a contradiction.

In other words, $\clos((\textsc{nor})\cap{\mathfrak{c}(\ttt{nil}_2)})\subsetneq(\textsc{nor})\cap\clos({\mathfrak{c}(\ttt{nil}_2)})$.
\end{eg}

%%%%%%%%%%%%%%%%%%%%%%%%%%%%%%%%%%%%%%%%%%

\subsection{An obstruction to belonging to $\ttt{clos}(\cniltwo)$} \label{sec4.17}
A complete characterisation of $\ttt{clos}(\cniltwo)$ requires us to understand which non-biquasitriangular operators lie in that set.   Although we do not have a complete answer to this question, we identify an index obstruction which allows us to eliminate large classes of operators from belonging to $\ttt{clos}(\cniltwo)$.

%%%%%%%%%%%%%%%%%%%%%%%%%%%%%%%%%%%%%%%%%%

\begin{prop} \label{prop4.18}
Let $T \in \bofh$ be a Fredholm operator with $\ttt{ind}(T)\in 2\mathbb{Z}+1$.   Then there exists $\delta > 0$ such that if $M \in \niltwo$ and $B \in \bofh$, then
\[
\norm T - (M B - B M) \norm \ge \delta. \]
In particular, $T \notin \ttt{clos}(\cniltwo)$.
\end{prop}

\begin{pf}
Since the set of Fredholm operators of index $\ttt{ind}(T)$ is open in $\bofh$, we can choose $\delta > 0$ such that if $C \in \bofh$ and $\norm C - T \norm < \delta$, then $\ttt{ind}(C) = \ttt{ind}(T) \in 2 \bbZ + 1$.

Suppose that  there exist $M \in \niltwo$ and $B \in \bofh$ such that $\norm T - (M B -  BM) \norm < \delta$.    Set $T_1 : = M B - B M$, implying that $\ttt{ind}(T_1) \in 2 \bbZ + 1$.

Writing $M = \begin{bmatrix} 0 & A \\ 0 & 0 \end{bmatrix}$ and $B =\begin{bmatrix} W & X \\ Y & Z \end{bmatrix}$ relative to $\ker\, M \oplus (\ker\, M)^\perp$ yields  $T_1 =\begin{bmatrix} AY &  AZ-WA\\ 0 & -YA \end{bmatrix}$. Since $T_1$ is Fredholm, $\dim \ker\, M =\dim  (\ker\, M)^\perp=\infty$, and $\pi(AY)$ is left invertible in the Calkin algebra
and $\pi(-YA)$ is right invertible. Hence $\pi(A)$ and $\pi(Y)$ are left and right invertible in the Calkin algebra, that is,
$\pi(A)$ and $\pi(Y)$ are invertible in the Calkin algebra. Therefore,
 $A$ and $Y$ are both Fredholm operators. By Corollary 5 of \cite{HanLeeLee2000},
\[\textrm{ind}(T_1)=\textrm{ind}(AY)+\textrm{ind}(-YA)=2[\textrm{ind}(A)+\textrm{ind}(Y)]\in 2\mathbb{Z},\]
a contradiction.

The second statement is an immediate consequence of the first.
\end{pf}

\begin{prop} \label{prop4.20}
Let $\{e_n\}_{n=1}^\infty$ be an orthonormal basis of $\cH$,  and let $W\in B(\cH)$ be the unilateral forward weighted shift with weight sequence $(w_n)_n \in \ell^\infty$.  Suppose that $\inf_n\, |w_n|>0$. Then
$W \notin \ttt{clos}(\cniltwo).$
\end{prop}

\begin{pf}
Without loss of generality, we may assume that $w_n>0$, $n\geq 0$.   Moreover,  $W= S D$, where $S$ is the  unilateral shift, and $D$ is the diagonal operator $\textup{diag} \{ w_n\}_n$.
Since $ \inf_n |w_n|>0$, $D$ is invertible. As $S$ is a Fredholm operator with index $-1$,  so is $W$. By Proposition~\ref{prop4.18}, $W\notin \ttt{clos}(\cniltwo).$
\end{pf}
\bigskip

%%%%%%%%%%%%%%%%%%%%%%%%%%%%%%%%%%%%%%%%%%
%%%%%%%%%%%%%%%%%%%%%%%%%%%%%%%%%%%%%%%%%%
%%%%%%%%%%%%%%%%%%%%%%%%%%%%%%%%%%%%%%%%%%
%%%%%%%%%%%%%%%%%%%%%%%%%%%%%%%%%%%%%%%%%%
%%%%%%%%%%%%%%%%%%%%%%%%%%%%%%%%%%%%%%%%%%
%%%%%%%%%%%%%%%%%%%%%%%%%%%%%%%%%%%%%%%%%%
%%%%%%%%%%%%%%%%%%%%%%%%%%%%%%%%%%%%%%%%%%
%%%%%%%%%%%%%%%%%%%%%%%%%%%%%%%%%%%%%%%%%%
%%%%%%%%%%%%%%%%%%%%%%%%%%%%%%%%%%%%%%%%%%
%%%%%%%%%%%%%%%%%%%%%%%%%%%%%%%%%%%%%%%%%%
%%%%%%%%%%%%%%%%%%%%%%%%%%%%%%%%%%%%%%%%%%
%%%%%%%%%%%%%%%%%%%%%%%%%%%%%%%%%%%%%%%%%%
%%%%%%%%%%%%%%%%%%%%%%%%%%%%%%%%%%%%%%%%%%
%%%%%%%%%%%%%%%%%%%%%%%%%%%%%%%%%%%%%%%%%%
%%%%%%%%%%%%%%%%%%%%%%%%%%%%%%%%%%%%%%%%%%
%%%%%%%%%%%%%%%%%%%%%%%%%%%%%%%%%%%%%%%%%%
%%%%%%%%%%%%%%%%%%%%%%%%%%%%%%%%%%%%%%%%%%
%%%%%%%%%%%%%%%%%%%%%%%%%%%%%%%%%%%%%%%%%%
%%%%%%%%%%%%%%%%%%%%%%%%%%%%%%%%%%%%%%%%%%
%%%%%%%%%%%%%%%%%%%%%%%%%%%%%%%%%%%%%%%%%%
%%%%%%%%%%%%%%%%%%%%%%%%%%%%%%%%%%%%%%%%%%
%%%%%%%%%%%%%%%%%%%%%%%%%%%%%%%%%%%%%%%%%%

\section{Concluding remarks.} \label{section5}

\subsection{} \label{sec5.01}
As mentioned at the outset, one may study the sets $\mathfrak{c}(\Omega) := \{ [A, B]: A, B \in \Omega\}$ and $\ttt{clos} (\mathfrak{c}(\Omega))$  for many classes $\Omega$ of operators on a Hilbert space.  In the infinite-dimensional setting, it is of interest to consider $\ttt{clos}(\mathfrak{c}(\Omega))$ whenever $\Omega$ is invariant under similarity, because in this case $\mathfrak{c}(\Omega)$ is also invariant under similarity, and we can make use of the machinery developed by Apostol, Fialkow, Herrero and Voiculescu~\cite{ApostolFialkowHerreroVoiculescu1984} to help us in this setting.

In the finite-dimensional setting, the existence of Jordan forms gives one hope that one may be able to describe $\mathfrak{c}(\Omega)$ itself when $\Omega$ is invariant under similarity.

%%%%%%%%%%%%%%%%%%%%%%%%%%%%%%%%%%%%%%%%%%

\subsection{} \label{sec5.02}
In this paper, we have considered the case where $\Omega = \niltwo$, which is clearly a similarity-invariant subset of $\bofh$, and in a previous paper~\cite{MarcouxRadjaviZhang2023}, we considered the case where $\Omega = \fE := \{E \in \bofh: E^2 = E\}$, again a similarity-invariant subset of $\bofh$.   As we saw in Theorem~\ref{thm4.11},
\[
\ttt{(bqt)} \cap \ttt{clos}({\cniltwo}) = \ttt{(bqt)} \cap \ttt{clos} (\fC_\fE). \]
Observe that square-zero operators and idempotents both satisfy a quadratic polynomial.   Suppose that
\[
\ttt{quad} := \{ T \in \bofh: q(T) = 0 \text{ for some polynomial } q  \text{ of degree at most } 2\}. \]
\begin{itemize}
	\item{}
	What are $\mathfrak{c}(\ttt{quad})$ and $\ttt{clos}(\mathfrak{c}(\ttt{quad}))$?
	\item{}
	%\textcolor{myred}{
Is $\ttt{(bqt)} \cap \ttt{clos}({\cniltwo}) = \ttt{(bqt)} \cap \ttt{clos}(\mathfrak{c}(\ttt{quad}))$?
\end{itemize}

%%%%%%%%%%%%%%%%%%%%%%%%%%%%%%%%%%%%%%%%%%

If we maintain the requirement that $\Omega \subseteq \bofh$ be invariant under similarity, but drop the requirement that elements of $\Omega$ must satisfy a polynomial of degree at most two, then -- as an extension of the work done here -- a natural class to consider is the set $\Omega = \ttt{nil} = \{ M \in \bofh:  M^k = 0 \text{ for some } k \ge 1\}$.   We have already seen that Smith's Theorem~\ref{thm3.18} implies that in the finite-dimensional setting, $\mathfrak{c}(\ttt{niil}) = \{ T \in \bbM_n(\bbC) : \ttt{Tr}(T) = 0\}$, and hence $\mathfrak{c}(\ttt{nil}) = \mathfrak{c}(\bbM_n(\bbC))$.

%%%%%%%%%%%%%%%%%%%%%%%%%%%%%%%%%%%%%%%%%%

When $\dim\, \hilb = \infty$, it is also true that the set $\ttt{clos}(\mathfrak{c}(\ttt{nil}))$ is strictly larger than $\ttt{clos}(\cniltwo)$.

\begin{prop} \label{prop5.03}
Every compact operator is a limit of commutators of compact, nilpotent operators.
\end{prop}

\begin{pf}
Let $K \in \kofh$ and $\eps > 0$.   Choose a finite-rank operator $F$ such that $\norm K - F \norm < \frac{\eps}{2}$, and write
\[
F= \begin{bmatrix} F_0 & 0 \\ 0 & 0 \end{bmatrix} \]
relative to $\hilb = \hilb_1 \oplus \hilb_2$, where $\dim\, \hilb_1 < \infty$.

Let $\alpha := \ttt{Tr}(F_0)$, and choose $m > 2$ such that $\frac{|\lambda|}{m} < \frac{\eps}{2}$.   Write $\hilb_2 := \hilb_{3} \oplus \hilb_4$, where $\dim \, \hilb_3 = m$, and -- relative to the decomposition $\hilb = \hilb_1 \oplus \hilb_3 \oplus \hilb_4$, define
\[
G := \begin{bmatrix} F_0 & 0 & 0 \\ 0 & \frac{-\lambda}{m} I_m & 0 \\ 0 & 0 & 0 \end{bmatrix}. \]
Then $\norm G- F \norm < \frac{\eps}{2}$ and therefore
\[
\norm K - G \norm \le \norm K - F \norm + \norm F - G \norm < \eps. \]
Since $F_0 \oplus \frac{-\lambda}{m} I_m \in \cB(\hilb_1 \oplus \hilb_3)$ has trace zero, it follows from Smith's Theorem~\ref{thm3.18} that $F_0 \oplus \frac{-\lambda}{m} I_m$ is a commutator of nilpotent matrices in $\cB(\hilb_1 \oplus \hilb_3)$, and thus $G$ is a commutator of finite-rank nilpotent operators.    But then $K$ is a limit of such operators.
\end{pf}

%%%%%%%%%%%%%%%%%%%%%%%%%%%%%%%%%%%%%%%%%%%

\subsection{} \label{sec5.04}
When $\Omega$ is not invariant under similarity, the problem of characterising $\mathfrak{c}(\Omega)$ and its closure appears to be more delicate.   Two interesting classes to consider are
\begin{enumerate}
	 \item[(a)]
	 $\Omega = \ttt{(nor)}$, the set of normal operators in $\bofh$, and
	 \item[(b)]
	 $\Omega = \cU$, the set of unitary operators in $\bofh$.
\end{enumerate}
We finish with two minor observations regarding this last class.

%%%%%%%%%%%%%%%%%%%%%%%%%%%%%%%%%%%%%%%%%%%

\begin{eg} \label{eg5.05}
Every $W \in \cU_2(\bbC)$ with $\ttt{Tr}(W) = 0$ lies in $\mathfrak{c}(\cU_2)$.

Note that if $W \in \cU_2(\bbC)$ and $\ttt{Tr}(W) = 0$, then $W$ is unitarily equivalent to $\begin{bmatrix} \theta & 0 \\ 0 & -\theta \end{bmatrix}$ and hence
\[
 W \simeq \theta \begin{bmatrix} 1 & 0 \\ 0 & -1 \end{bmatrix} \simeq  \theta \begin{bmatrix} 0 & 1 \\ 1 & 0 \end{bmatrix} \]
for some $\theta \in \bbT$.

Let $U_0 = \begin{bmatrix} 1 & 0 \\ 0 & i \end{bmatrix}$, $V_0 = \frac{(1+i)}{2} \begin{bmatrix} 1 & 1 \\ -1 & 1 \end{bmatrix}$, so that $U_0$ and $V_0$ are unitary.  Then
\[
W_0 := [U_0, V_0]  = \frac{(1+i)}{2} \begin{bmatrix} 0 & 1-i \\ 1-i & 0 \end{bmatrix} = \begin{bmatrix} 0 & 1 \\ 1 & 0 \end{bmatrix}.\]

Let $U = \theta U_0$ and $V = V_0$.  It follows that $W \simeq [U, V]$ is a commutator of unitary operators.
\end{eg}

%%%%%%%%%%%%%%%%%%%%%%%%%%%%%%%%%%%%%%%%%%

\begin{eg} \label{eg5.06}
There does not exist a unitary operator in $\cB(\bbC^3)$ which is a commutator of unitary operators in $\cB(\bbC^3)$.

\bigskip

Clearly it suffices to show that if $U = \begin{bmatrix} \alpha_1 & 0 & 0 \\ 0 & \alpha_2 & 0 \\ 0 & 0 & \alpha_3 \end{bmatrix}$ with respect to some orthonormal basis $\{ e_1, e_2, e_3\}$ for $\bbC^3$ and $V = [v_{ij}]$ are unitary operators in $\cB(\bbC^3)$, then $W := [U, V]$ is not unitary.

Suppose otherwise.

If $\alpha_i = \alpha_j$ for some $i \ne j$, then $(U- \alpha_i I)$ has rank at most one, and so
\[
W = [U,V] = [U - \alpha_i I, V] \]
has rank at most two, contradicting the fact that $W$ is unitary.  Thus $\alpha_1, \alpha_2, \alpha_3$ are distinct.

The hypothesis that $W = [w_{ij}] = [(\alpha_i - \alpha_j) v_{ij}]$ is unitary implies that its columns form an orthonormal set.   Consider
\[
0 = \left \langle \begin{bmatrix} 0 \\ (\alpha_2 - \alpha_1) v_{21} \\ (\alpha_3 - \alpha_1) v_{31} \end{bmatrix}, \begin{bmatrix} (\alpha_1 - \alpha_2) v_{12} \\ 0 \\ (\alpha_3 - \alpha_2) v_{32} \end{bmatrix} \right \rangle = (\alpha_3 - \alpha_1) v_{31} \ol{(\alpha_3 - \alpha_2) v_{32}}. \]
Since the $\alpha_i$'s are all distinct, this implies that $v_{31} = 0$ or $v_{32} = 0$.  By interchanging the first two basis vectors $e_1$ and $e_2$ if necessary (which doesn't affect the fact that $U$ is diagonal), we may assume without loss of generality that $v_{31} = 0$.

Next, consider the first and third columns of $W$;
\[
0 = \left \langle \begin{bmatrix} 0 \\ (\alpha_2 - \alpha_1) v_{21} \\ 0  \end{bmatrix}, \begin{bmatrix} (\alpha_1 - \alpha_3) v_{13} \\ (\alpha_2-\alpha_3) v_{23} \\ 0 \end{bmatrix} \right \rangle = (\alpha_2 - \alpha_1) v_{21} \ol{(\alpha_2 - \alpha_3) v_{23}}. \]
Arguing as before, since the $\alpha_i$'s are all distinct, this implies that $v_{21} = 0$ or $v_{23} = 0$.

\begin{itemize}
	\item{}
	If $v_{21} = 0$, then the first column of $W$ is zero, contradicting the hypothesis that $W$ is unitary.
	\item{}
	If $v_{23} = 0$, then $V = \begin{bmatrix} v_{11} & v_{12} & v_{13} \\ v_{21} & v_{22} & 0 \\ 0 & v_{32} & v_{33} \end{bmatrix}$.
	But $V$ is unitary, and so the last two rows of $V$ are mutually orthogonal.   From this we deduce that $v_{22} = 0$ or $v_{32} = 0$.   If $v_{32} = 0$, then a quick calculation shows that the last row of $W$ is zero, again contradicting the hypothesis that $W$ is unitary.   Thus $v_{22} = 0$.   Since $V$ is unitary, its first and second rows are orthogonal, proving that $v_{11}=0$.
	
	 The hypothesis that $W$ is unitary then implies that its first and third rows are orthogonal, which in turn implies that $v_{12} = 0$ or $v_{32}= 0$.   From above, this in turn forces $v_{12} = 0$.   Since $V$ is also unitary, this then forces $v_{33} = 0$.  We conclude that
	 \[
	 V = \begin{bmatrix} 0 & 0 & v_{13} \\ v_{21} & 0 & 0 \\ 0 & v_{32} & 0 \end{bmatrix}. \]
\end{itemize}	
In particular, $|v_{13}| = |v_{21}| = |v_{32}| = 1$.

Hence
\[
W =  \begin{bmatrix} 0 & 0 & (\alpha_1 - \alpha_3) v_{13} \\ (\alpha_2 - \alpha_1) v_{21} & 0 & 0 \\ 0 & (\alpha_3 - \alpha_2) v_{32} & 0 \end{bmatrix}. \]
Since $W$ is unitary, we again conclude that $|(\alpha_1 - \alpha_3) v_{13}| = | (\alpha_2 - \alpha_1) v_{21}| = | (\alpha_3 - \alpha_2) v_{32}| = 1$, and so $|(\alpha_1 - \alpha_3)| = | (\alpha_2 - \alpha_1)| = | (\alpha_3 - \alpha_2) |= 1$.

In other words, $\alpha_1, \alpha_2$ and $\alpha_3 \in \bbT$ are all at a distance of one from one another.   Clearly this is impossible.

\end{eg}

%%%%%%%%%%%%%%%%%%%%%%%%%%%%%%%%%%%%%%%%%%

\begin{eg} \label{eg5.07}
The above phenomenon does not extend to higher dimensions.    Let $W = \mathrm{diag} (\alpha, -\alpha, \beta, -\beta) \in \cU_4(\bbC)$.   Then $W = W_1 \oplus W_2$ where, by Example~\ref{eg5.05}, $W_1 = [U_1, V_1]$ and $W_2 = [U_2, V_2]$ for unitary matrices $U_1, U_2, V_1, V_2 \in \cU_2(\bbC)$.  Hence
\[
W = [U, V ], \]
where $U = U_1 \oplus U_2$ and $V= V_1 \oplus V_2$.
\end{eg}

%%%%%%%%%%%%%%%%%%%%%%%%%%%%%%%%%%%%%%%%%%

%%%%%%%%%%%%%%%%%%%%%%%%%%%%%%%%%%%%%%%%%%
%%%%%%%%%%%%%%%%%%%%%%%%%%%%%%%%%%%%%%%%%%
%%%%%%%%%%%%%%%%%%%%%%%%%%%%%%%%%%%%%%%%%%
%%%%%%%%%%%%%%%%%%%%%%%%%%%%%%%%%%%%%%%%%%
%%%%%%%%%%%%%%%%%%%%%%%%%%%%%%%%%%%%%%%%%%
%%%%%%%%%%%%%%%%%%%%%%%%%%%%%%%%%%%%%%%%%%
%%%%%%%%%%%%%%%%%%%%%%%%%%%%%%%%%%%%%%%%%%

%%%%%%%%%%%%%%%%%%%%%
% Bibliography
%%%%%%%%%%%%%%%%%%%%%

\bibliographystyle{plain}
%\bibliography{/Users/mjdruitt/Library/texmf/bibtex/2020papers}

\end{document}